\documentclass[12pt,reqno]{amsart}

\input{amssym.def}
\input{amssym}

\def\be{\begin{equation}}
\def\ee{\end{equation}}
\def\ba{\begin{eqnarray}}
\def\ea{\end{eqnarray}}
\def\lb{\label}
\def\nn{\nonumber}
\def\ni{\noindent}
\def\cal{\mathcal}

\def\D{{ D_{\hspace{-2pt}\raisebox{-2pt}{\tiny\sl R}}}}
\def\Dop{{ D_{\hspace{-2pt}\raisebox{-2pt}{$\scriptstyle R_{op}$}}}}
\def\Rpsi{{\Psi_{\hspace{-2pt}\raisebox{-2pt}{\tiny\sl R}}}}
\def\tr{{\rm Tr}\,}
\def\Rtr{{\rm Tr}_{\hspace{-2pt} \raisebox{-2pt}{\tiny\sl R}\,}}
\def\Rinvtr{{\rm Tr}_{\hspace{-2pt} \raisebox{-2pt}{$\scriptstyle R^{-1}$}\,}}
\def\Roptr{{\rm Tr}_{\hspace{-2pt} \raisebox{-2pt}{$\scriptstyle R_{op}$}\,}}
\def\RTr#1{{\rm Tr}_{\hspace{-2pt} \raisebox{-2pt}{\tiny\sl R}\,\raisebox{2pt}{\scriptsize$(#1)$}}}
\def\RopTr#1{{\rm Tr}_{\hspace{-2pt} \raisebox{-2pt}{$\scriptstyle R_{op}$}\,\raisebox{2pt}{\scriptsize$(#1)$}}}
\def\Tr#1{{\rm Tr}_{\raisebox{2pt}{\scriptsize$(#1)$}}}
\def\Rdet{{{\rm det}_{\hspace{-1pt} \raisebox{-2pt}{\tiny\sl R}\,}}}

\makeatletter
\@addtoreset{equation}{section}\makeatother

\newcounter{theorem}
\makeatletter\@addtoreset{theorem}{section}\makeatother

\newtheorem{prop}[theorem]{Proposition}
\newtheorem{rem}[theorem]{Remark}

\newtheorem{def-lem}[theorem]{Definition-Lemma}
\newtheorem{def-prop}[theorem]{Definition-Proposition}
\newtheorem{theor-def}[theorem]{Theorem-Definition}
\newtheorem{defin}[theorem]{Definition}
\newtheorem{theor}[theorem]{Theorem}
\newtheorem{cor}[theorem]{Corollary}

\begin{document}
\begin{titlepage}
\title[Unitary quantum group differential calculus]
{On construction of unitary quantum group differential calculus}

\author{Pavel Pyatov}

\address{Pavel Pyatov
Laboratory of Mathematical Physics,
National Research University Higher School of Economics
20 Myasnitskaya street, Moscow 101000, Russia \&
Bogoliubov Laboratory of Theoretical Physics, Joint Institute
for Nuclear Research, 141980 Dubna, Moscow Region, Russia}
\email{pyatov@theor.jinr.ru}

\thanks{
This research was supported by the National Research
University –- `Higher School of Economics’ Academic Fund Program (grant No.14-01-0027 for the period 2014-2015) and by the
grant of RFBR No.14-01-00474-a.}

\date{\today}

\begin{abstract}
We develop a construction of the unitary type anti-involution for the quantized differential calculus over $GL_q(n)$ in the case $|q|=1$. To this end, we consider a joint associative algebra of quantized functions, differential forms and Lie derivatives over $GL_q(n)/SL_q(n)$, which is bicovariant 
with respect to $GL_q(n)/SL_q(n)$ coactions. We define a specific non-central {\em spectral extension} of this algebra by the spectral variables of 
three matrices of the algebra generators. In the spectrally expended algebra we construct a three-parametric family of its inner automorphisms.
These automorphisms are used for construction of the unitary anti-involution for the (spectrally extended) calculus over $GL_q(n)$. 
\end{abstract}

\maketitle

\end{titlepage}

\phantom{a}
\vspace{2cm}
\tableofcontents
\vfill\newpage

\section{Introduction}
\lb{sec1}

Soon after the invention of the quantum group theory \cite{Drin} the construction of the differential calculi over quantum spaces and quantum groups
became a hot topic in the noncommutative geometry accumulating much investigation activity. The general frameworks
for these investigations were given by the bicovariance postulates \cite{Wor} and the R-matrix ideology \cite{FRT}. A substantial progress
was soon achieved in constructing the algebras of differential operators (Lie derivatives) and differential forms over the general linear quantum groups
(see, e.g., \cite{Jur,Malt,Sud,Tzy,SWZ.92,SWZ.93,Z}). At the same time, serious difficulties were met in all attempts of complementing the calculus
with more sophisticated structures. This concerns studies of the quantum group's real forms, where no unitary quantum groups were found for the most interesting regime of the quantization parameter  $|q|=1$. This happened to a construstion of the exterior algebra for the differential forms over quantum orthogonal and symplectic groups, where the no-go theorem was proved (see \cite{AIP}). This was also the case with the construction of de Rahm complex over special linear quantum groups in the frame of the Woronowicz's approach. 
So it comes out that in the quantum group calculus introducing any additional structure meets serious troubles and  demands case-by-case investigations. 

In the present paper, we address the problem of construction of the unitary anti-involution for the differential calculus over linear quantum groups.
One crucial hint for its solution was given already in  \cite{AF.92} where the  anti-involution map was looked for and found not in the quantum group, but in the larger Heisenberg double algebra.   Another  important ingredient of our construction is the spectral extension of the calculus algebra
which was elaborated in \cite{IP.09} on the basis of the Cayley-Hamilton theorem for the quantum matrix algebras (see \cite{GPS.97,IOP.98,IOP.99}).

The paper is organized as follows. 
In the next section we describe the differential calculi algebras over general and special linear quantum groups.
We introduce algebras of quantized functions, differential forms and Lie derivatives over $GL_q(n)$ and consider their $SL_q(n)$ reduction and their bicovariant structure. We follow mainly the papers \cite{SWZ.92,IP.95}, though we consider several different sets of generators and disuss the $SL_q(n)$ reduction conditions in detail.

In section \ref{sec3}, we recall briefly the Cayley-Hamilton theorem for the Reflection equation algebras of the $GL_q(n)$ type. These are the algebras of the left- and right-invariant Lie derivatives. Extending the ideas of \cite{IP.09} we construct a special non-central extension of the differential calculus algebra with spectral variables  --- the eigenvalues of  three matrices generating the Reflection equation subalgebras in the differential calculus (two of them are the matrices of the left- and right-invariant Lie derivatives). This is one of the two main results of this paper, we present it in theorem \ref{theor-def3.2}. Then, in the spectrally extended algebra we  introduce  a three-parametric family of inner automorphisms. For certain integer values of their parameters these automorphisms reproduce discrete time evolution of the q-deformed top \cite{AF.92}. 

Section \ref{sec4} describes the construcion of Gauss decomposition for the Reflection equation algebras. 
Starting from this section we restrict our consideration to the algebras associated with the Drinfeld-Jimbo R-matrices (see (\ref{DJ}) below).
All the previous constructions were carried out for a more general family of R-matrices of the $SL(n)$ type. Explanation of
this notion and a collection of the R-matrix formulas are given in Appendix.

Section \ref{sec5} is devoted to construction of the unitary anti-involution in the (spectrally extended) differential calculus over $GL_q(n)$.
This is the second major result of this work and we present it in theorem \ref{theor5.5}. Restriction of the anti-involution to the  differential calculus over $SL_q(n)$ is also discussed.  

We finally note that although the exterior derivative is not considered in this paper, its construction suggested in \cite{FP.94} looks appropriate for the calculi we discuss. We leave a detailed consideration of the exterior derivative, its compatibility with the unitary structure and its possible BRST realization for a future publication.

\section{Differential forms and Lie derivatives on $GL_q(n)$}
\lb{sec2}

In this section we describe associative unital algebras of differential calculi over $GL_q(n)$ and $SL_q(n)$
---  quantizations of the classical calculi over $GL(n)$/$SL(n)$. The quantum calculi algebras are generated by the components of four $n\times n$ matrices:
\begin{itemize}
\item[]
${\cal k}T^i_j{\cal k}_{i,j=1}^n$ -- coordinate functions,\smallskip
\item[]
${\cal k}\Omega^i_j{\cal k}_{i,j=1}^n$ -- right-invariant 1-forms,\smallskip
\item[]
${\cal k}L^i_j{\cal k}_{i,j=1}^n$ -- right-invariant Lie derivatives,\smallskip
\item[]
${\cal k}K^i_j{\cal k}_{i,j=1}^n$ -- left-invariant Lie derivatives.
\end{itemize}
In the $GL_q(n)$ case one imposes on the generators a set of quadratic relations which fix classical values\footnote{E.g., ${n^2 +k \choose k}$ and, respectively, ${n^2\choose k}$ are dimensions of the spaces of $k$-th order homogeneous polynomials of coordinate functions and, respectively, right-invariant 1-forms.} for the dimensions of the spaces of homogeneous polynomials in generators.
These relations, in general, allow alphabetic ordering of the generators and therefore,
we call them {\em permutation relations}. Transition to the $SL_q(n)$ calculus is then achieved by imposing one more polynomial
relation for each matrix of generators. We call these relations {\em reduction conditions}.

Both, the permutation relations and the reduction conditions are given with the use of R-matrix ---
an $n^2\times n^2$ matrix $R\in {\rm Aut}({\Bbb C}^n\otimes {\Bbb C}^n)$ satisfying Artin's braid relation.
We specify $R$ to be of the {\it $SL(n)$ type} which means, in particular, that
its minimal polynomial is quadratic.
All the necessary notions and the basics of the R-matrix techniques are recalled in the Appendix.
For a more detailed presentation of the R-matrix formalism the reader is referred to paper \cite{IP.09} and references therein.
In what follows we use notation adopted there.

Our main motivating example is the calculus associated with the standard Drinfeld-Jimbo R-matrix
\be
\lb{DJ}
R\, =\, \sum_{i,j=1}^n q^{\delta_{ij}}E_{ij}\otimes E_{ji}\, +\,
(q-q^{-1})\,\sum_{i<j } E_{ii}\otimes E_{jj}\, .
\ee
Here $q\in {\Bbb C}\setminus\{0\}$,~ and
$(E_{ij})_{kl}:=\delta_{ik}\delta_{jl}$, $\, i,j=1,\dots ,n,\,$ are matrix units.
As we will show  the quantum differential calculus associated with the matrix (\ref{DJ})
admits unitary type specialization. This is the major result of our work.
However, we stress that consistent $SL_q(n)$ and $GL_q(n)$ differential calculi
can be associated with any $SL(n)$ type R-matrix.
Among those are multiparametric generalizations of the Drinfeld-Jimbo R-matrix \cite{R.90} (see also example 2.10 in \cite{IP.09}) and the
Cremmer-Gervais R-matrices \cite{CG.90,H}.
Therefore, we present a part of the construction in full generality and stick to considering the particular Drinfeld-Jimbo's case
starting from section \ref{sec4}.

We now proceed to writing explicit formulas.

\subsection{Quantization of functions}\hfill \par\vspace{3pt}
\lb{subsec2.1}

\noindent
Throughout this and the next sections we assume that $R\in {\rm Aut}({\Bbb C}^n\otimes {\Bbb C}^n)$ is the $SL(n)$ type R-matrix.
This notion together with the compressed matrix index notation and the notions of {\it R-trace},  ~$\Rtr$,~ and of
{\it q-antisymmet\-rizers}, ~$A^{(n)}$,~ which appear in formulas below are explained in the Appendix.

Permutation relations for the quantized coordinate functions over $GL_q(n)$ in the compressed matrix index notation read \cite{Drin,FRT}
\be
\lb{RTT}
R \, T_1 \, T_2 = T_1 \, T_2 \, R\, .
\ee
Strictly speaking, one also assumes invertibility of {\it quantum determinant} of the matrix $T\,$
\be
\lb{detT}
\Rdet T := \Tr{1, \dots ,n} \bigl( A^{(n)} T_1 T_2 \dots T_n \bigr)\, .
\ee
The $SL(n)$ type  property of the R-matrix $R$ guarantees centrality of $\Rdet T$
and the $SL_q(n)$ reduction condition reads
\be
\lb{sl-1}
\Rdet T\, =\, 1.
\ee
Further on we denote an associative unital algebra generated by elements $T_{ij}$ subject to relations (\ref{RTT}), (\ref{sl-1})
as ${\cal F}[R]$ and call it {\em the algebra of functions over $SL_q(n)$}. It is also briefly called {\em the RTT algebra}.
\smallskip

In a standard way ${\cal F}[R]$ is endowed with the Hopf algebra structure \cite{Drin,FRT} from which we recall formulas for
the coproduct and for the antipode:
\ba
\lb{bi-RTT}
\Delta (T_{ij}) &=& \textstyle \sum_{k=1}^n T_{ik}\otimes T_{kj}\, ,
\\[2pt]
\lb{antipode}
(T^{-1})_1 &=&
q^{n(n-1)}\, n_q\, \RTr{2,\dots ,n}\left(T_2\dots T_n A^{(n)}\right) (\Rdet T)^{-1}\, ,
\ea\smallskip
where $n_q:= (q^n-q^{-n})/( q-q^{-1})$ is the q-number.
Using the symbol $T^{-1}$ for the antipode instead of the standard notation $S(T)$ is justified by equalities\smallskip
\be
\lb{T-inverse}
\textstyle
\sum_{k=1}^n T_{ik}\, (T^{-1})_{kj}\, =\, \sum_{k=1}^n (T^{-1})_{ik}\, T_{kj}\, =\,
\delta_{ij}\, 1\, .
\ee
\smallskip

\subsection{Quantization of forms}\hfill \par\vspace{3pt}
\lb{subsec2.2}

\noindent
Permutation relations for the quantized external algebra of the right-invariant differential forms over $GL_q(n)$ 
were first suggested in \cite{SWZ.92,Z} 
\be
\lb{OmegaOmega-GL}
R\,\Omega^g_1 R\, \Omega^g_1\, =\, -\, \Omega^g_1 R\, \Omega^g_1 R^{-1}\, .
\ee
This algebra implies following permutation relations for the R-trace of $\Omega^g$
$$
\left(\Rtr\Omega^g\right) \Omega^g\, +\, \Omega^g\left(\Rtr\Omega^g\right)\, =\, -(q-q^{-1})\, (\Omega^g)^2\,
$$
thus making a naive $SL_q(n)$ reduction $\Rtr \Omega^g =0$ impossible. Instead, one observes
that the R-traceless matrix
\be
\lb{O-SL}
\Omega\, :=\, \Omega^g - \frac{q^n}{n_q} \left(\Rtr\Omega^g\right) I
\ee
generates a subalgebra in (\ref{OmegaOmega-GL}) which does not contain element $\Rtr \Omega^g$:
\be
\lb{OmegaOmega}
R\,\Omega_1 R\, \Omega_1 \, +\, \Omega_1 R\, \Omega_1 R^{-1} \, =\,
\kappa_q \left( \Omega_1^2\, +\, R\,\Omega_1^2 R\right)\, .
\ee
Here we denote
\be
\lb{kappa}
\kappa_q :=\frac{ q^n (q-q^{-1})}{ n_q + q^n  (q-q^{-1}) },\qquad\mbox{assuming additionally}\qquad
n_q + q^n  (q-q^{-1})\neq 0.
\ee

Permutation relations (\ref{RTT}),(\ref{OmegaOmega}) complemented with
\be
\lb{OmegaT}
R\,\Omega_1 R^{-1} T_1 \, =\, T_1\, \Omega_2 \,
\ee
define the algebra which is consistent with the $SL_q(n)$ type reduction conditions  for functions (\ref{sl-1})
and for 1-forms
\be
\lb{Omega-sl}
\Rtr \Omega\, =\, 0\, .
\ee
We  call it {\em the external algebra of differential forms over} $SL_q(n)$ \cite{IP.95}.
Analogously, relations (\ref{RTT}),(\ref{OmegaOmega-GL}) complemented with (\ref{OmegaT}), where $\Omega$ is to be substituted by $\Omega^g$,
define {\em the external algebra of differential forms over} $GL_q(n)$. As we will see, all relations containing matrices of 1-forms
$\Omega^g/\Omega$ in the differential calculi over $GL_q(n)/SL_q(n)$ are identical, with the only exception of their own permutations 
(\ref{OmegaOmega-GL}) and (\ref{OmegaOmega}). Later on we will always write down relations for $\Omega$ understanding that 
their analogues for $\Omega^g$ have identical form.
\smallskip

Equivalently, the external algebra can be written in terms of the {\em left-invariant} 1-forms
\be
\lb{Theta}
\Theta^g_{ij}\,  =\, 
(T^{-1}\Omega^g T)_{ij}, \qquad
\Theta_{ij}\,  =\, (T^{-1}\Omega T)_{ij},
\ee
Deriving the permutation relations and the reduction condition for $\Theta$ is a good exercise in the R-matrix techniques.
They read:
\ba
\lb{ThetaTheta-GL}
&&
R^{-1}\Theta^g_2 R\, \Theta^g_2 \, =\, -\, \Theta^g_2 R\, \Theta^g_2 R  ,
\\[2pt]
\lb{ThetaTheta}
&&
R^{-1}\Theta_2 R\, \Theta_2 \, +\, \Theta_2 R\, \Theta_2 R \, =\, \kappa_q \left( \Theta_2^2\, +\, R\,\Theta_2^2 R\right) ,
\\
\lb{Theta-sl}
&&
{\rm Tr}_{\hspace{-2pt} \raisebox{-2pt}{\tiny\sl $R_{o p}$\!\!}}\Theta \, =\, 0 ,
\\[2pt]
\lb{ThetaT}
&&
\Theta_1 T_2  \, =\, T_2 R^{-1}\Theta_2 R , \qquad \mbox{\em (for both, $\Theta$ and $\Theta^g$).}
\ea
Here $R_{op} := P R P$, ~and $P\in {\rm Aut}({\Bbb C}^n\otimes {\Bbb C}^n)$ is the permutation operator: $P(u\otimes v) = v\otimes u$.
\smallskip

{\em Bi-invariant} objects in the external algebra are given by the R-traces of powers of $\Omega$.
Their subalgebra is not affected by the quantization
\cite{FP.96} and looks like
\ba
\lb{TrO-even}
\Rtr \Omega^{2i}& =& {\rm Tr}_{\hspace{-2pt} \raisebox{-2pt}{\tiny\sl $R_{o p}$\!\!}}\Theta^{2i}\, =\,0 \qquad\qquad\; \forall\, i\geq 1,
\\[2pt ]
\lb{TrO-odd}
\omega_i \, :=\, \Rtr \Omega^{2i+1} & =&
{\rm Tr}_{\hspace{-2pt} \raisebox{-2pt}{\tiny\sl $R_{o p}$\!\!}}\Theta^{2i+1}\, \;\;\,\forall\, i=1,\dots ,n-1,
\\[2pt]
\omega_i\,  \omega_j& =& -\, \omega_j\,  \omega_i\, .
\ea

\subsection{Quantization of Lie derivatives}\hfill \par\vspace{3pt}
\lb{subsec2.3}

\noindent
The algebra of quantized right-invariant Lie derivatives (see eq.(\ref{REA}) below) is widely known under the
name of {\it Reflection equation (RE)} algebra. It was introduced in a context of a factorized particle scattering on a half-line \cite{Ch,KS} and since,
has found a number of applications in the theory of integrable systems, in noncommutative geometry and in quantum groups.

The common algebra of quantized functions and invariant Lie derivatives given by eq.(\ref{RTT}) and by eqs.(\ref{REA}),(\ref{LT}) below has also its own name ---
the {\it Heisenberg double (HD)} algebra \cite{AF.91,STS,SWZ.93}. The action of Lie derivatives on functions is identical to the action
of the corresponding vector fields and so, the HD algebra can be considered as the algebra of differential operators over quantum group.

The whole set of relations for the right-invariant Lie derivatives was written for the first time in \cite{SWZ.92,Z}:
\ba
\lb{REA}
R\, L_1 R\, L_1 & =& L_1 R\, L_1 R \, ,
\\
\lb{LT}
R\, L_1 R\, T_1 & =&  \gamma^{2} \, T_1 L_2\, ,
\\
\lb{LOmega}
R\, L_1 R\, \Omega_1 & = &  \Omega_1 R\, L_1 R\, .
\ea

For the $SL_q(n)$ reduction we introduce quantum determinant of $L$ as\footnote{Note that the difference in 
the definitions of the quantum determinants for the
matrices $T$ and $L$ comes from the difference in their permutation relations.
}
\be
\lb{detL}
\Rdet L \, := \, \RTr{1, \dots ,n} \left( A^{(n)}
L_{\overline{1}} \dots L_{\overline{n}} \right) ,
\ee
where concise notation
\be
\lb{L-over}
L_{\overline{1}}:= L_1\, , \quad L_{\overline{i}} := R_{i-1} \, L_{\overline{i-1}} \, R_{i-1}^{-1}\quad \forall\; i> 1.
\ee
for the matrix $L$ $\overline{i}$-th copy is used. 
Further on, we
assume invertibility and demand centrality of $\Rdet L$ --- the latter condition fixes value of the parameter  $\gamma$ in (\ref{LT})
\be
\lb{gamma-SL}
\gamma = q^{1/n}.
\ee
The $SL_q(n)$ reduction condition then reads:
\be
\lb{L-sl}
\Rdet L\, = \, q^{-1}\, 1.
\ee
A particular convenience of the normalization factor $q^{-1}$ chosen here will become obvious later on (see section \ref{subsec3.3}, the definition of the automorphism $\phi_L$).

\begin{rem}
\lb{rem-Lie-GL}
In the algebra of Lie derivatives over $GL_q(n)$ the element $\Rdet L$ can not be central. Usually this is achieved by the choise
$\gamma =1$. Alternatively, one can keep $\gamma = q^{2/n}$ extending the algebra with one more invertible generator $\ell$, satisfying
permutation relations
$$
\ell\, L\, =\, L\, \ell, \qquad \ell\, T\, =\, q^{2/n}\, T\, \ell .
$$
Later on, we will not consider separately the algebra of Lie derivatives over $GL_q(n)$ having this possibility in mind.
\end{rem}

Let's now turn to discussion of the left-invariant Lie derivatives.
Unlike vector fields the left- and right-invariant Lie derivatives
are independent and have to be introduced in the calculus separately.
Notice however that, keeping the mirror (left-right) symmetry of the calculus
one can uniquely reproduce permutation relations for the
left-invariant objects ($K$, $\Theta$, etc.) reading relations for right-invariant generators literally leftwards
(note that the mirror image of $R=R_{12}$ is $R_{op}=R_{21}$). In this way for the left-invariant Lie derivatives $K$ we obtain
\ba
\lb{REA-K}
R\, K_2 R\, K_2 & =& K_2 R\, K_2 R \, ,
\\
\lb{KT}
T_2\, R\, K_2 R & =&  \gamma^{2} \, K_1 T_2\, ,
\\
\lb{ThetaK}
R\, K_2 R\, \Theta_2 & =& \Theta_2\, R\, K_2 R \, .
\ea
These relations are to be complemented by the natural commutativity conditions of
the left/right-invariant Lie derivatives with all the right/left-invariant objects:
\ba
\lb{K-X}
K_1\, X_2\, = \,  X_2 \, K_1\, && \forall\, X\, =\, \Omega, L ,
\\
\lb{L-Y}
L_1\, Y_2 \, = \,  Y_2 \, L_1\,  && \forall\,  Y\, =\, \Theta, K  .
\ea
Here the only new relation is the commutativity condition for $L$ and $K$.
All the other commutativity conditions follow from the permutation relations imposed earlier.

The $SL_q(n)$ reduction for the left-invariant Lie derivatives $K$ reads
\be
\lb{detK-sl}
\Rdet K\, := \,
{\rm Tr}_{\hspace{-2pt} \raisebox{-2pt}{\tiny\sl $R_{op}$}\,\raisebox{2pt}{\scriptsize$(1, \dots ,n)$}}
\left(A^{(n)}
K_{\underline{n}} \dots K_{\underline{1}} \right)
\, = \, q^{-1}\, 1,
\ee
where $\underline i$-th copy of the matrix $K$ is defined as
\be
\lb{K-under}
K_{\underline{n}}:= K_n\, , \quad K_{\underline{i}} := R_i \, K_{\underline{i+1}} \, R_i^{-1}\quad \forall\; 1\leq i< n.
\ee
Here again, formulae for $\Rdet K$ are obtained as mirror copies of those for $\Rdet L$.

\subsection{Yet another RE algebra}\hfill \par\vspace{3pt}
\lb{subsec2.5}

\noindent
Ideologically, permutation relations for the left- and right-invariant Lie derivatives
given in the previous section describe an effect of, respectively,
right and left shifts on the underlying quantum group `manifold'.
It would be instructive to study distinctions in their actions on (i.e., permutations with) the other objects of the calculus.
To this end one may compare two objects obeying the same invariance properties: $L$ and $TK^{-1}T^{-1}$.
It turns out that they have identical permutation relations with functions, but not with forms.
Let us analyze quantitatively this phenomenon looking at permutation relations of their ratio
\be
\lb{F}
F_{ij}\, :=\, (L\, T\, K\, T^{-1})_{ij}\, .
\ee

\begin{prop}
\lb{F-algebra}
The components of matrix ${\cal k}F^i_j{\cal k}_{i,j=1}^n$ generate the RE subalgebra in the calculi algebra
and satisfy following permutation relation with the other generators
\ba
\lb{REA-F}
R\, F_1 R\, F_1 & =& F_1 R\, F_1 R \, ,
\\
\lb{FT}
R\, F_1 R^{-1} T_1 & =&   T_1 F_2\, ,
\\
\lb{FOmega}
R\, F_1 R\, \Omega_1 & = &  \Omega_1 R\, F_1 R\, ,
\\
\lb{FL}
R\, L_1 R\, F_1 & =& F_1 R\, L_1 R \, ,
\\
\lb{FK}
F_1\, K_2 & = &  K_2 \, F_1\, .
\ea
The $SL_q(n)$ reduction (\ref{gamma-SL}),(\ref{L-sl}),(\ref{detK-sl}) implies following condition on $F$
\be
\lb{det-F}
\Rdet F \, := \, \RTr{1, \dots ,n} \left( A^{(n)}
F_{\overline{1}} \dots F_{\overline{n}} \right)\, =\, q^{-n^2} 1 .
\ee
\end{prop}

\ni
{\bf Proof.}
Checking permutation relations (\ref{REA-F})-(\ref{FK}) is straightforward and we skip it.
We shall consider in details a more sophisticated calculation of the quantum determinant of $F$.

First of all, denoting $U=T K T^{-1}$ we separate $L$ and $U$ factors in $\Rdet F$:
$$
\Rdet F\, =\, \RTr{1, \dots ,n}  \left(A^{(n)}\left( L_{\overline{1}} \dots L_{\overline{n}} \right)
\left( U_{\overline{1}} \dots U_{\overline{n}} \right)\right) .
$$
Here since matrix $U$ satisfies a version of the reflection equation with the inverse matrix $R$ in it:
$R^{-1} U_1 R^{-1} U_1 =  U_1 R^{-1} U_1 R^{-1}$,~
the definition of its overlined copy  differs from that
for $L$: $U_{\overline{i}}:= R_{i-1}^{-1} U_{\overline{i-1}} R_{i-1}$ (c.f. with eq.(\ref{L-over})).
Using condition ${\rm rk} A^{(n)}=1$ (see Appendix) we simplify the expression
$$
\Rdet F\, =\,
\Rdet L \cdot
\Tr{1, \dots ,n} \bigl( A^{(n)}  U_{\overline{1}} \dots U_{\overline{n}}  \bigr) .
$$
Next, using permutation relations for $T$ and $K$ we separate different matrix factors in the product of $U$-copies:
$$
U_{\overline{1}} \dots U_{\overline{n}} \, =\, {\gamma}^{-2n(n-1)}\left({\textstyle \prod_{i=1}^nJ_i}\right)^2
\left(T_1\dots T_n\right) \left(K_{\underline{n}}\dots K_{\underline{1}}\right) \left(T^{-1}_n\dots T^{-1}_1\right).
$$
Here
$J_i$ are the R-matrix realizations of a commutative set of Jucys-Murphy elements
in the braid group (see, e.g., \cite{OP.01}):
\be
\lb{Jucys}
J_1\, :=\, I, \qquad J_{i+1}\, :=\, R_i J_i R_i\quad\forall\, i\geq 1.
\ee

Evaluating $(\prod_{i=1}^nJ_i)$ at $A^{(n)}$ as $q^{-n(n-1)}$ and exploiting again the  rank=1 property of $A^{(n)}$
we obtain
$$
\Rdet F\, =\, (\gamma q)^{-2n(n-1)}\,
\Rdet L \cdot  \Rdet T \cdot \Tr{1, \dots ,n} \bigl( A^{(n)}  K_{\underline{n}}\dots K_{\underline{1}}  \bigr) \cdot  (\Rdet T)^{-1} .
$$
Finally, using eq.(\ref{n3}) we substitute ${\rm Tr}$ in this formula by ${\rm Tr}_{\hspace{-2pt} \raisebox{-2pt}{\tiny\sl $R_{op}$}}$
and get
\ba
\nn
\Rdet F& =&\gamma^{-2n(n-1)}\, q^{-n(n-2)}\,\Rdet L \cdot  \Rdet T\cdot \Rdet K \cdot  (\Rdet T)^{-1}
\\[2pt]
\nn
&=& q^{-n^2}\, \gamma^{2n}\,\Rdet L \cdot \Rdet K\, =\, q^{-n^2} 1.
\ea
Here in the last line we use relation
$$
T\cdot \Rdet K \, =\, q^{-2}\, \gamma^{2n}\, \Rdet K \cdot T
$$
to cancel $\Rdet T$  and then, apply
the $SL_q(n)$ reduction conditions.
\hfill$\blacksquare$
\smallskip

\subsection{Bicovariance}\hfill \par\vspace{3pt}
\lb{subsec2.4}

\noindent
Up to now we have described certain unital associative algebra generated by the components of four matrices $T$, $\Omega$, $L$ and $K$. A remarkable fact which makes this algebra indeed the differential calculus over the quantum group is a possibility to endow it with a
structure of the bicovariant bimodule\footnote{For definition of the bicovariant bimodule see, e.g., \cite{KSch}} over Hopf algebra ${\cal F}[R]$.
In this section we complement the construction of the differential calculus describing the ${\cal F}[R]$ comodule structures over it.

The left and right ${\cal F}[R]$ coactions --- ~$\delta_\ell$ and $\delta_r$~ --- on the algebra generators are defined as follows:

\noindent
- on $T_{ij}$ they  just reproduce the coproduct (\ref{bi-RTT})
\be
\lb{coact-rl-T}
\delta_{\ell/r} (T_{ij}) \, =\, \Delta (T_{ij}) ;
\ee

\noindent
- on the matrices of right-invariant generators, such as $\Omega$, $L$, or $F$, the coactions are given by
\be
\lb{coact-r-inv}
\textstyle
\delta_\ell (X_{ij})\, =\, \sum_{k,p=1}^n(T_{ik}\otimes 1)(1\otimes X_{kp}) ((T^{-1})_{pj}\otimes 1) , \qquad
\delta_r (X_{ij})\, =\, X_{ij}\otimes 1 , \quad
\ee
where  $X\, =\, \Omega,\, L,\, F,\, \dots$;
\medskip

\noindent
- on the matrices of left-invariant generators, such as $\Theta$, or $K$, they are defined as
\be
\lb{coact-l-inv}
\textstyle
\delta_\ell (Y_{ij})\, =\, 1\otimes Y_{ij}, \qquad \delta_r(Y_{ij})\, = \, \sum_{k,p=1}^n(1\otimes (T^{-1})_{ik})(Y_{kp}\otimes 1) (1\otimes T_{pj}),
\ee
where  $Y\, =\, \Theta,\, K,\,  \dots$.
\medskip

\noindent
The use of terminology "left/right-invariant" becomes now evident. \smallskip

Notice that the co-transformation
properties of the generators are preserved under matrix multiplication (e.g., $\sum_k L_{ik} L_{kj}$ and $\sum_k \Theta_{ik} K_{kj}$ are right- and left-invariant, respectively),
whereas conjugation with $T$ interchanges left and right co-transformations (e.g., $\sum_{k,p}T_{ik} K_{kp} T^{-1}_{pj}$ and $\sum_{k,p} T^{-1}_{ik} L_{kp} T_{pj}$
transform as right- and left-invariant objects, respectively).
The operation of taking R-trace extracts bi-invariant objects, so that for $\omega_i$ from eq.(\ref{TrO-odd}) one has
$$
\delta_\ell (\omega_i)\, =\, 1\otimes \omega_i,\qquad \delta_r (\omega_i)\, =\, \omega_i\otimes 1 .
$$

\subsection{Summary}\hfill \par\vspace{3pt}
\lb{subsec2.6}

\noindent
We collect considerations of this section into a
\begin{defin}
\lb{defin-algebra}
To any $SL(n)$-type R-matrix $R$ there corresponds an
associative unital algebra  ${\goth D\goth C}_{gl}[R]$ of the differential calculus over $GL_q(n)$.
This algebra is generated by the components of four matrices
$T$, $\Omega^g$, $L$ and $K$,
subject to the permutation relations
\begin{center}
(\ref{RTT}), (\ref{OmegaOmega-GL}), (\ref{OmegaT}), (\ref{REA})-(\ref{LOmega}), (\ref{REA-K}), (\ref{KT}), (\ref{K-X}),
\end{center}
Substituting in this definition generators $\Omega^g\mapsto \Omega$  and relations (\ref{OmegaOmega-GL})$\mapsto$(\ref{OmegaOmega})    
and adding the $SL_q(n)$ reduction conditions
\begin{center}
(\ref{sl-1}), (\ref{Omega-sl}), (\ref{kappa}), (\ref{gamma-SL}), (\ref{L-sl}), (\ref{detK-sl}).
\end{center}
one obtains definition of the algebra  ${\goth D\goth C}_{sl}[R]$ of the differential calculus over $SL_q(n)$.
The bicovariant ${\cal F}[R]$-bimodule structure on both  algebras is given by
eqs.(\ref{coact-rl-T})-(\ref{coact-l-inv}).
\end{defin}

\begin{rem}
As a generating set for ${\goth D\goth C_{gl/sl}}[R]$ one can choose also quadruples of
matrices $\{T, \Theta^g/\Theta, L, K\}$ and $\{T, \Omega^g/\Omega, L, F\}$. The permutation relations and the reduction conditions for these sets were
presented earlier in this section.
\end{rem}

\begin{rem}
The Heisenberg double algebra over $SL_q(n)$ investigated in \cite{IP.09} is a quotient algebra of ${\goth D\goth C}_{sl}[R]$
over relations
\be
\lb{red-HD}
F_{ij}\, =\, \delta_{ij}\, 1, \qquad \Omega_{ij}\, =\, 0.
\ee
The first relation imposes dependence of the left- and right-invariant vector fields in the Heisenberg double.
\end{rem}

\section{Spectral extension and automorphisms}
\lb{sec3}

In this section we introduce three families of automorphisms on algebras  ${\goth D\goth C}_{gl/sl}[R]$.
Two of these automorphisms are generated by the actions of Lie derivatives $L$, $K$ and, as explained in \cite{AF.92}, they
reproduce a q-deformed version of an evolution of the Euler's isotropic top.
The third automorphism is related with matrix $F$ and it acts on forms leaving functions invariant.
In section 5 we use these automorphisms for
construction of the unitary anti-involution over ${\goth D\goth C}_{gl}[R]$ and so, we have to define them as the algebra inner automorphisms.
To this end we define a {\em spectral extension} of the algebra --- its extension by the eigenvalues of matrices $L$, $K$ and $F$.
For the Heisenberg double algebra the spectral extension was constructed in \cite{IP.09}. Here we present generalization of that construction
to the algebras ${\goth D\goth C}_{gl/sl}[R]$.
\smallskip

\subsection{Characteristic identities and spectral variables}\hfill  \par\vspace{3pt}
\lb{subsec3.1}

\noindent
In this subsection we collect structure results about the RE algebras of the $SL(n)$ type \cite{GPS.97,IOP.98,IOP.99}.
These results are necessary for the subsequent constructions.

Consider the RE algebra (\ref{REA}), (\ref{L-sl}) generated by the matrix of right-invariant Lie derivatives $L$. A set of elements
$a_i$, $i=0,\dots, n$,
\be
\lb{aL}
a_0\, :=\, 1\, , \qquad
a_i \, := \, \RTr{1, \dots ,i} \left( A^{(i)}
L_{\overline{1}} \dots L_{\overline{i}} \right)\, , \quad i\geq 1.
\ee
belongs to the center of the RE algebra; the last of them -- $a_n$ -- is just the quantum determinant of $L$.
These elements are the coefficients of the following matrix identity
\be
\lb{charL}
\sum^{n}_{i=0}\, (-q)^i\, a_{i}\, L^{n-i}\, =\, 0 \, ,
\ee
which is nothing but the RE algebra analogue of the  Cayley-Hamilton theorem.
We will introduce a special central extension of the RE algebra with the aim at bringing
the {\em characteristic identity}  (\ref{charL}) to a factorized form.

Consider an Abelian $\Bbb C$-algebra of polynomials in $n$ invertible
indeterminates $\{\mu_\alpha^{\pm 1}\}_{\alpha =1}^n$ and in their differences
$\{(\mu_\alpha-\mu_\beta)^{\pm 1}\}_{\alpha >\beta=1}^n$, satisfying
condition
$$
\textstyle
\prod_{\alpha=1}^n \mu_\alpha = q^{-1}.
$$
We parameterize elements $a_i$ of the RE algebra by the elementary symmetric polynomials in $\mu_\alpha$:
\be
\lb{aL-elementary}
a_i = e_i(\mu_1,\dots ,\mu_n):= \sum_{1\leq \alpha_1<\dots <\alpha_i\leq n}
\mu_{\alpha_1} \mu_{\alpha_2}\dots \mu_{\alpha_i}
\qquad
\forall\, i=0,1,\dots ,n\, ,
\ee
assuming commutativity of indeterminates $\mu_\alpha$ with the elements of the RE algebra
\be
\lb{L-mu}
L\, \mu_\alpha\, =\, \mu_\alpha\, L\, .
\ee
The resulting central extension of the RE algebra is called its {\it spectral extension}, and the elements $\mu_\alpha$
are called {\it eigenvalues of the `quantum' matrix} $L$.
The characteristic identity in the completed RE algebra assumes a factorized form
\be
\lb{charL-factor}
\prod_{\alpha =1}^n\bigl(L - q\mu_{\alpha} I\bigr)\, =\, 0\, .
\ee
It can be used for the construction of a set of mutually orthogonal matrix idempotents

\be
\lb{projL}
P^\alpha\, :=\, \prod_{\beta=1\atop \beta\neq\alpha}^n\frac{\bigl(L - q\mu_{\beta} I\bigr)}{ q(\mu_\alpha-\mu_\beta)}\; :\quad
P^\alpha P^\beta\, =\, \delta_{\alpha\beta}\, P^\alpha\, , \quad
\sum_{\alpha=1}^n P^\alpha\, =\, I\, .
\ee
By construction, evaluating  $L$ on the idempotents one obtains the eigenvalues
\be
\lb{L-P}
L\, P^\alpha\, =\, P^\alpha L\, =\, q\mu_\alpha P^\alpha\, .
\ee

Now we apply similar spectral extension procedure for the RE algebras generated by matrices $K$ and $F$: see, respectively, eqs.(\ref{REA-K}), (\ref{detK-sl}) and (\ref{REA-F}), (\ref{det-F}).
We parameterize coefficients $b_i$ and $c_i$ of their characteristic identities by
elementary symmetric functions in $n$ indeterminates $\nu_\alpha$ and $\rho_\alpha$,  $\alpha=1,\dots ,n$:
\ba
\lb{bK}
b_0\, :=\, 1\, , &&
b_i \, := \,{\rm Tr}_{\hspace{-2pt} \raisebox{-2pt}{\tiny\sl $R_{op}$}\,\raisebox{2pt}{\scriptsize$(1, \dots ,i)$}}
\left( A^{(i)} K_{\underline{i}} \dots K_{\underline{1}} \right) \, =\, e_i(\nu_1,\dots ,\nu_n)\, ,
\\
\lb{cF}
c_0\, :=\, 1\, , &&
c_i \, := \, \RTr{1, \dots ,i} \left( A^{(i)}
F_{\overline{1}} \dots F_{\overline{i}} \right)
\, =\, e_i(\rho_1,\dots ,\rho_n)\, ,
\ea
where
$$
\textstyle
\prod_{\alpha=1}^n \nu_\alpha = q^{-1}, \qquad \prod_{\alpha=1}^n \rho_\alpha = q^{-n^2}.
$$
Assuming centrality of the eigenvalues
\be
\lb{KF-nu-rho}
\nu_\alpha\, K\, =\, K\, \nu_\alpha\, ,\qquad \rho_\alpha\, F\, =\, F\, \rho_\alpha\, ,
\ee
we factorize the characteristic identities
\ba
\lb{charK}
\sum^{n}_{i=0}\, (-q)^i\, b_{i}\, K^{n-i} &=& \prod_{\alpha =1}^n\bigl(K - q\nu_{\alpha} I\bigr)\, =\, 0 ,
\\
\lb{charF}
\sum^{n}_{i=0}\, (-q)^i\, c_{i}\, F^{n-i} &=& \prod_{\alpha =1}^n\bigl(F - q\rho_{\alpha} I\bigr)\, =\, 0 .
\ea
Imposing additionally invertibility conditions on the eigenvalues and on their differences,
we define associated sets of matrix idempotents
\ba
\lb{projK}
Q^\alpha &:=& \prod_{\beta=1\atop \beta\neq\alpha}^n\frac{\bigl(K - q\nu_{\beta} I\bigr)}{ q(\nu_\alpha-\nu_\beta)}\; :\quad
Q^\alpha Q^\beta\, =\, \delta_{\alpha\beta}\, Q^\alpha\, , \quad
\sum_{\alpha=1}^n Q^\alpha\, =\, I\, ,
\\
\lb{projF}
S^\alpha &:=& \prod_{\beta=1\atop \beta\neq\alpha}^n\frac{\bigl(F - q\rho_{\beta} I\bigr)}{ q(\rho_\alpha-\rho_\beta)}\; :\quad
S^\alpha S^\beta\, =\, \delta_{\alpha\beta}\, S^\alpha\, , \quad
\sum_{\alpha=1}^n S^\alpha\, =\, I\, ,
\ea
so that
\ba
\lb{K-Q}
K\, Q^\alpha &=& Q^\alpha K\, =\, q\nu_\alpha Q^\alpha\, ,
\\[2pt]
\lb{F-S}
F\, S^\alpha &=& S^\alpha F\,\, =\, q\rho_\alpha S^\alpha\, .
\ea

Finally, we can consistently set that the newly introduced spectral variables
transform trivially under both left anf right ${\cal F}[R]$ coactions
\be
\lb{coact-munurho}
\delta_r \xi \, =\, \xi \otimes 1 , \;\;
\delta_\ell \xi \, =\, 1\otimes \xi ,\qquad \forall\, \xi\in \{\mu_\alpha,\nu_\alpha,\rho_\alpha\} .
\ee
\smallskip

\subsection{Spectral extension 
}\hfill  \par\vspace{3pt}
\lb{subsec3.2}

\noindent
Our next step is to  construct  an extension of  algebras
${\goth D\goth C}_{gl/sl}[R]$ with the spectral vaiables. By no means the result is going to be a trivial central extension.
Our goal is to define permutation relations for $\mu_\alpha$, $\nu_\alpha$ and $\rho_\alpha$ in such a way,
that their elementary symmetric functions
would commute with $T$ and $\Omega$ exactly as the elements $a_i$, $b_i$ and $c_i$ do.
In \cite{IP.09} a consistent definition for permutations of $\mu_\alpha$ with $T$ was derived. Here we follow the same scheme.
First, we calculate permutation relations of $a_i$, $b_i$ and $c_i$ with $T$ and $\Omega$.

\begin{prop}
\lb{prop3.1}
In algebras ${\goth D\goth C}_{gl/sl}[R]$ elements $a_i$ (\ref{aL}) satisfy permutation relations
\ba
\lb{Ta_i}
\gamma^{2i}\, T \, a_i  &=&
\textstyle
a_i \, T\, -\,
(q^2-1) \sum_{j=1}^{i} (-q)^{-j} a_{i-j}\, (L^j  T)\, ,
\\[5pt]
\lb{Oa_i}
\bigl[\Omega\, ,\, a_i\bigr] &=&
\textstyle
(q^2-1)\, \sum_{j=1}^i (-q)^{-j}\bigl[\Omega\, , a_{i-j} L^j \bigr]\, .
\ea
Here notation $[\,\cdot\, ,\,\cdot\, ]$ stands for the commutator.
Relations for elements $b_i$ (\ref{bK}) are mirror images of those for $a_i$ with the substitution
$a_i\mapsto b_i$, $L\mapsto K$, $\Omega\mapsto \Theta$.
Namely,
\ba
\lb{Tb_i}
\gamma^{2i}\,  b_i\, T   &=&
\textstyle
T\, b_i \,-\,
(q^2-1) \sum_{j=1}^{i} (-q)^{-j}  ( T K^j)\, b_{i-j}\,,
\\[3pt]
\lb{Ob_i}
\bigl[\Theta\, ,\, b_i\bigr] &=&
\textstyle
(q^2-1)\, \sum_{j=1}^i (-q)^{-j}\bigl[\Theta\, , b_{i-j} K^j \bigr]\, .
\ea
Permutation relations of $c_i$ (\ref{cF}) with $\Omega$ are identical to (\ref{Oa_i}) with the substitution
$a_i\mapsto c_i$, $L\mapsto F$, while with $T$ elements $c_i$ commute.
\end{prop}

\ni
{\bf Proof.} Eq.(\ref{Ta_i}) is proved in \cite{IP.09} in Proposition 3.18.
The proof of eq.(\ref{Oa_i}) is based on equality
\ba
\lb{fromCHN}
q^{i(i-1)}\RTr{2,\dots ,i+1}\Bigl((L_{\overline{1}}J_1)\dots
(L_{\overline{i}} J_i)\,  A^{(i)}\Bigr)^{\uparrow 1}\!\! =
a_i I_1 - (q^2-1) \sum_{j=1}^i (-q)^{-j} a_{i-j} L_1^j , \quad
\ea
where $X^{\uparrow 1} := I\otimes X\in {\rm End}(V^{\otimes (i+1)})\;\; \forall\, X\in {\rm End}(V^{\otimes i})$,
and the Jucys-Murphy elements $J_i$ were defined in (\ref{Jucys}).

Equality (\ref{fromCHN}), in turn, follows from the Cayley-Hamilton-Newton identity (see Theorem 3.11 in \cite{IP.09}) and it is
contained implicitly in the proof of Proposition 3.18  in \cite{IP.09}.

Now, as a consequence of permutation relations (\ref{LOmega}) one has
$$
(L_{\overline j} J_j) \Omega_1\, =\, \Omega_1 (L_{\overline j} J_j), \quad \forall\, j\geq 2,
$$
and, hence the l.h.s of (\ref{fromCHN}) commutes with $\Omega_1$. So does the r.h.s., which immediately leads to the equality (\ref{Oa_i}).

Permutation relations for $b_i$  (\ref{Tb_i}), (\ref{Ob_i}) follow from   (\ref{Ta_i}), (\ref{Oa_i})
and the left-right symmetry of the calculus.

Identical form of the commutators of $a_i$ and $c_i$ with $\Omega$ is a consequence of the identity of the
pemutation relations for $L$ and $F$ with $\Omega$. Checking relation $[c_i, T]=0$ is straightforward.
\hfill$\blacksquare$

\begin{theor-def}
\lb{theor-def3.2}
Consider an Abelian $\Bbb C$-algebra of polynomials in $3n$  invertible
indeterminates and in their differences
$$
\{\mu_\alpha^{\pm 1},\,\nu_\alpha^{\pm 1},\,\rho_\alpha^{\pm 1},\,(\mu_\alpha-\mu_\beta)^{\pm 1},\,(\nu_\alpha-\nu_\beta)^{\pm 1},\,
(\rho_\alpha-\rho_\beta)^{\pm 1}\}_{\alpha >\beta=1}^n,
$$
subject to relations
$$
\textstyle
\prod_{\alpha=1}^n \mu_\alpha\, = \,\prod_{\alpha=1}^n \nu_\alpha\, =\, q^{-1}, \qquad \prod_{\alpha=1}^n \rho_\alpha\, =\, q^{-n^2}.
$$

Spectral extensions $\widehat{\goth D\goth C}_{gl/sl}[R]$ of the $GL_q(n)/SL_q(n)$ differential calculi ${\goth D\goth C}_{gl/sl}[R]$
by this algebra
are given by parameterization formulas  (\ref{aL-elementary}), (\ref{bK}), (\ref{cF}) and  by permutation relations
\ba
\lb{X-mu}
\mu_\alpha\, X\, =\, X\, \mu_\alpha\, , &&  \forall\, X= L,\, K,\, \Theta\, , G:= T^{-1} F T,
\\[2pt]
\lb{Y-nu}
\nu_\alpha\, Y\, =\, Y\, \nu_\alpha\, , &&  \forall\, Y= L,K,\Omega\, ,F ,
\\[2pt]
\lb{Z-rho}
\rho_\alpha\, Z\, =\, Z\, \rho_\alpha\, , &&  \forall\, Z= T,L,K,F,G ,
\ea
\vspace{-8mm}
\ba
\lb{T-mu}
\gamma^2 \, (P^\beta T)\, \mu_\alpha & = &
q^{2\delta_{\alpha\beta}} \mu_\alpha\, (P^\beta T) ,
\\[2pt]
\lb{T-nu}
\gamma^2 \, \nu_\alpha\, (T Q^\beta) & = &
q^{2\delta_{\alpha\beta}}  (T Q^\beta )\, \nu_\alpha, \quad \mbox{\small (recall that $\gamma=q^{1/n}$),}
\\[2pt]
\lb{Omega-mu}
 q^{2\delta_{\alpha\sigma}} (P^\beta\, X\, P^\sigma)\, \mu_\alpha &=&
q^{2\delta_{\alpha\beta}}  \mu_\alpha\, (P^\beta\, X\, P^\sigma)\, \quad \forall X=\Omega, F,
\\[3pt]
\lb{Theta-nu}
 q^{2\delta_{\alpha\sigma}} (Q^\beta\, Y\, Q^\sigma)\, \nu_\alpha &=&
q^{2\delta_{\alpha\beta}}  \nu_\alpha\, (Q^\beta\, Y\, Q^\sigma)\, \quad\, \forall Y=\Theta, G,
\\[3pt]
\lb{Omega-rho}
 q^{2\delta_{\alpha\sigma}} (S^\beta\, \Omega\, S^\sigma)\, \rho_\alpha &=&
q^{2\delta_{\alpha\beta}}  \rho_\alpha\, (S^\beta\, \Omega\, S^\sigma)  \quad\;\; \forall\,
\alpha ,\beta , \sigma =1,\dots ,n\, .
\\
\lb{weaker}
\left[\rho_\alpha, \, \Rtr X\right]\, =\, 0, && \mbox{\it where $X$ is any matrix monomial in $\Omega$ and $F$.}
\ea
Here expressions for matrix idempotents $P^\alpha$, $Q^\alpha$, $S^\alpha$ are given in (\ref{projL}), (\ref{projK}), (\ref{projF}).

Formulae (\ref{T-mu})-(\ref{Omega-rho}) can be equivalently written as
\ba
\lb{T-mu2}
\gamma^2 \,  T\, \mu_\alpha & =&
\mu_\alpha\, T\, +\, (q-q^{-1})\, (L P^\alpha T)\, ,
\\[2pt]
\lb{T-nu2}
\gamma^2 \, \nu_\alpha\,  T  & =&
 T\, \nu_\alpha\, +\, (q-q^{-1})\, (T K Q^\alpha)\, ,
\\[2pt]
\lb{Omega-mu2}
\left[ X\, ,\, \mu_\alpha \right] &=& (q-q^{-1})\, \left[ L P^\alpha\, ,\, X\right]\, \quad\;\,\, \forall X=\Omega,F,
\\[2pt]
\lb{Theta-nu2}
\left[Y\, ,\, \nu_\alpha \right] &=& (q-q^{-1})\, \left[ K Q^\alpha\, ,\, Y\right]\, \quad\;\, \forall Y=\Theta,G,
\\[2pt]
\lb{Omega-rho2}
\left[ \Omega\, ,\, \rho_\alpha \right] &=& (q-q^{-1})\, \left[ F S^\alpha\, ,\, \Omega\right]\qquad \forall\alpha  =1,\dots ,n.
\ea

The extended algebras  $\widehat{\goth D\goth C}_{gl/sl}[R]$
are endowed with the structure of ${\cal F}[R]$ bicovariant bimodule by eqs.(\ref{coact-rl-T})-(\ref{coact-l-inv}), (\ref{coact-munurho}).
\end{theor-def}

\ni
{\bf Proof.}
As concerns formulae (\ref{T-mu})-(\ref{Omega-rho}),
one has to check that they are consistent with the parameterization
formulas and the relations obtained in  Proposition \ref{prop3.1}.
Relations (\ref{Ta_i}) were checked for consistency in \cite{IP.09}, Theorem 3.27.
Here we shall prove consistency of  the spectral extension with relations (\ref{Oa_i}).
The rest of relations follow by similar considerations.

Let us calculate permutation of the elementary symmetric function in spectral values $e_i(\mu)$ with the matrix of 1-forms $\Omega$ with the use of eq.(\ref{Omega-mu}). In calculations below we use following notations
$~~e_i(\mu^{/\alpha}):=e_i(\mu)|_{\mu_\alpha =0},\;\;\; e_i(\mu^{/\alpha \beta}):=e_i(\mu)|_{\mu_\alpha =\mu_\beta =0}$.
\ba
\nonumber
&&e_i(\mu)\, \Omega \, =\, \sum_{\alpha,\beta =1}^n e_i(\mu) P^\alpha \Omega P^\beta \, =\,\sum_\alpha P^\alpha \Omega P^\alpha e_i(\mu)
\\
\nonumber
&+&
\sum_{\alpha\neq \beta} P^\alpha \Omega P^\beta \Bigl(
e_i(\mu^{/\alpha \beta})+e_{i-1}(\mu^{/\alpha \beta})(q^{-2}\mu_\alpha+q^2\mu_\beta)+e_{i-2}(\mu^{/\alpha \beta})\mu_\alpha\mu_\beta\Bigr)
\ea
\ba
\nonumber
&=& \Omega\,  e_i(\mu)\, +\, (q^2-1)\sum_{\alpha\neq \beta}\left(
P^\alpha \Omega P^\beta\, e_{i-1}(\mu^{/\alpha \beta}) \mu_\beta\, -\,  e_{i-1}(\mu^{/\alpha \beta})\mu_\alpha\, P^\alpha \Omega P^\beta\right)
\\
\nonumber
&=& \Omega\,  e_i(\mu)\, +\,(q^2-1)\sum_{\alpha,\beta =1}^n\left(P^\alpha \Omega P^\beta\, e_{i-1}(\mu^{/\beta}) \mu_\beta\, -\,
 e_{i-1}(\mu^{/\alpha}) \mu_\alpha\, P^\alpha \Omega P^\beta\right)
\\
\nonumber
&=& \Omega\,  e_i(\mu)\, +\,(q^2-1)\sum_{\alpha =1}^n\left[ \Omega\, ,\, P^\alpha\, e_{i-1}(\mu^{/\alpha}) \mu_\alpha\right]
\\
\nonumber
&=& \Omega\,  e_i(\mu)\, +\,(q^2-1)\sum_{j=1}^{i}\left[ \Omega\, ,\,e_{i-j}(\mu)(-L/q)^j  \right].
\ea
Here in the last line one uses formula $e_i(\mu^{/\alpha})=\sum_{j=0}^i e_{i-j}(\mu)(-\mu_\alpha)^j$ and takes into account that $\mu_\alpha$ in the presence of $P^\alpha$ can be substituted by $L/q$.

Comparing the first and the last lines of the calculation we see consistency of the spectral extension with eq.(\ref{Oa_i}).

Finally, it is easy to see that eq.(\ref{weaker}) agrees with the algebraic relations in ${\goth D\goth C}[R]$
observing that, by eqs.(\ref{REA-F}), (\ref{FOmega}), all components of matrix $F$ commute with elements $\Rtr X$ from (\ref{weaker}).
\hfill$\blacksquare$
\smallskip

Strictly speaking, relations (\ref{T-mu})-(\ref{Omega-rho}) for spectral variables are not the pemutations, since they are non-quadratic.
Quite remarkably,  they are consistent with the  commutators (\ref{X-mu})-(\ref{Z-rho}).
For instance, the permutation relations of $\mu_\alpha$ with $T$ (\ref{T-mu2}) and  with $\Omega$ (\ref{Omega-mu2}) leed to the
trivial commutator for $\mu_\alpha$ and $\Theta = T^{-1} \Omega T$.

\begin{rem}
\lb{rem3.3}
Spectral variables $\mu_\alpha$, $\nu_\alpha$  satisfy stronger version of equality (\ref{weaker}): they commute with all bi-invariant elements of the calculus. In particular,
\be
\lb{mu-biinv}
\left[\xi, \, \Rtr X\right]\, =\, \bigl[ \xi, \, \Roptr \! Y\bigr]\, =\, 0\quad \forall \xi\in\{\mu_\alpha,\nu_\alpha\},
\ee
where $X/Y$ could be any matrix monomial in $\{\Omega,L,F\}/\{\Theta,K,G\}$. This relations follow from the commutativity (\ref{X-mu}),
(\ref{Y-nu}) and the fact that the  R-trace of any monomial in right-invariant matrices $L$, $\Omega$. $F$ can be reexpressed in terms of
$R_{op}$-traces of left-invariant matrices $K$, $\Theta$, $G$, and vice-versa.
\end{rem}

\begin{rem}
\lb{rem3.4}
Eq.(\ref{weaker}) in the definitions of $\widehat{\goth D\goth C}_{gl/sl}[R]$
can be equally substituted by condition
\be
\lb{opca}
\Rtr \bigl(S^\alpha X S^\beta\bigr) \, =\, 0\quad \forall\beta\neq \alpha,
\ee
where $X$ is any matrix monimial in $\Omega$ and $F$.
Indeed:
\ba
\nn
&
q\rho_\alpha \Rtr \bigl(S^\alpha X S^\beta\bigr)\, =\, \Rtr \bigl(F S^\alpha X S^\beta\bigr) \, =\,
\RTr{1,2} \bigl( F_1 R\, S^\alpha_1  X_1 S^\beta_1\bigr)
&
\\
\nn
&
\RTr{1,2} \bigl(R^{-1} S_1^\alpha X_1 S_1^\beta R F_1 R\bigr)\, =\, \Rtr \bigl( S^\alpha X S^\beta F\bigr)\, =\, \Rtr
\bigl(S^\alpha X S^\beta\bigr) q \rho_\beta,
&
\ea
wherefrom (\ref{opca}) follows,  if one commutes $\rho_\beta$ to the left and uses invertibility of  $(\rho_\beta - \rho_\alpha)$.
The opposite implication follows from the presentation $\Rtr X=\sum_\beta \Rtr (S^\beta X S^\beta)$ and the permutation relations
(\ref{Omega-rho}).
\end{rem}

\subsection{Automorphisms}\hfill  \par\vspace{3pt}
\lb{subsec3.3}

\noindent
In \cite{AF.92} an important discrete sequence of the Heisenberg double algebra automorphisms was introduced. This sequence generated by the right-invariant Lie derivatives was interpreted there as a discrete time evolution of the q-deformed Euler top and so, a problem of a construction of its evolution operator was posed. The problem was further addressed in \cite{IP.09}, where it was shown that a solution can be found after the spectral extension of the initial algebra. Moreover, in a extended algebra one has a continuous one-parametric family of the automorphisms --- a continuous time evolution.

In the differential calculi algebras ${\goth D\goth C}_{gl/sl}[R]$ one can define three independent series of such type automorphisms.
\begin{prop}
\lb{prop3.5}
Mappings $\phi_L$, $\phi_K$ and $\phi_F$, defined on  generators as\footnote{Normalizations  (\ref{L-sl}), (\ref{detK-sl}) of the matrices $L$ and $K$ were chosen in such way that the transformation rules for $T$ here do not contain nontrivial coefficients.}
\ba
\lb{phiL}
\phi_L \!\!\! &:& T \mapsto LT  , \;\;\;  \Omega \mapsto L\Omega L^{-1}, \;\;\;   F \mapsto L F L^{-1} , \;\;
X \mapsto X  \;\; \forall X = L, K, \Theta  ;
\\[2pt]
\lb{phiK}
\phi_K \!\!\! &:& T \mapsto TK  , \;\; \Theta \mapsto K^{-1}\Theta K   , \; Y \mapsto Y\;\; \forall Y = L, K, \Omega ,F  ;
\\[2pt]
\lb{phiF}
\phi_F \!\!\! &:& T \mapsto T  , \;\;\;\;\; \Omega \mapsto F\Omega F^{-1}   , \;\; Z \mapsto Z\;\; \forall Z = L, K, F  .
\ea
generate the algebra automorphisms of the differential calculi ${\goth D\goth C}_{gl/sl}[R]$. These automorphisms are mutually commutative.
\end{prop}

\ni
{\bf Proof.}
Checking compliance of the maps with the permutation relations in  ${\goth D\goth C}_{gl/sl}[R]~$ and their mutual commutativity is straightforward.
In Proposition 4.1 \cite{IP.09} mappings $\phi_L$ and $\phi_K$ are proved to comply with the $SL_q(n)$ the reduction conditions on the Lie derivatives. It lasts to test transformations of the reduction condition for the differential forms. We consider a calculation for $\phi_F(\Rtr\Omega)$:
\ba
\nn
\phi_F(\Rtr\Omega) &=& \Rtr\left( F \Omega F^{-1}\right) = \RTr{1,2}\left(R^{-1} \underline{R F R} \Omega F^{-1}\right) =
\RTr{1,2} \left( R^{-1} \Omega F^{-1} \underline{R F R}\right)
\\[2pt]
\nn
&=&
\RTr{1,2}\left( \Omega F^{-1} R F\right) = \Rtr\left( \Omega F^{-1} F\right) = \Rtr\Omega = 0,
\ea
where the underlined expression in the first line is moved to the right with the use of permutation relations for $F$.
\hfill$\blacksquare$
\smallskip

In the spectrally completed algebras $\widehat{\goth D\goth C}_{gl/sl}[R]~$ these mappings
can be generalized to a three-parametric family of inner algebra automorphisms.
Strictly speaking,  to this end one has to extend further the calculus, passing from
spectral generators $\{\mu_\alpha,\nu_\alpha,\rho_\alpha\}$
to a new set of variables $\{x_{\alpha},y_\alpha, z_\alpha\}_{\alpha=1}^n$
\be
\lb{log}
\mu_\alpha \, =\, q^{-1/n} \exp (2\pi i x_\alpha), \quad \nu_\alpha \, =\, q^{-1/n} \exp(2\pi i y_\alpha), \quad
\rho_\alpha \, =\, q^{-n} \exp(2\pi i z_\alpha),
\ee
and considering formal power series in $x_\alpha$, $y_\alpha$, $z_\alpha$.
In terms of these new variables the $SL_q(n)$ reduction conditions read
\be
\lb{sl-xy}
\textstyle
\sum_{\alpha=1}^{n} x_\alpha\, =\, \sum_{\alpha=1}^{n} y_\alpha\, =\, \sum_{\alpha=1}^{n} z_\alpha\, =\, 0,
\ee
and the permutation relations  (\ref{T-mu})-(\ref{Theta-nu})
take an additive form. For instance, permutations of $x_\alpha$ with the matrices $T$ and $\Omega$ read
\ba
\lb{T-x}
(P^\beta T)\, x_\alpha &=& \left(x_\alpha +2\tau(\delta_{\alpha\beta}-n^{-1})\right) (P^\beta T),
\\[2pt]
\lb{Omega-x}
(P^\beta \Omega P^\sigma)\, x_\alpha &=&
\left(x_\alpha +2\tau \delta_{\alpha\beta}-2\tau \delta_{\alpha\sigma}\right)(P^\beta \Omega P^\sigma),
\ea
where we denote
\be
\lb{tau}
\tau\, :=\, \frac{1}{ 2\pi i} \log q\, .
\ee
\smallskip

The main result of this section is the following
\begin{theor}
\lb{theor3.6}
Consider three-parametric family $\phi_{(t_1,t_2,t_3)}$ of  $\widehat{\goth D\goth C}_{gl/sl}[R]$  inner automorphisms
\ba
\nonumber
\phi_{(t_1,t_2,t_3)}&:& u\mapsto \varphi_{(t_1,t_2,t_3)}\, u\, (\varphi_{(t_1,t_2,t_3)})^{-1}\quad \forall\, u\in \widehat{\goth D\goth C}[R],
\\[4pt]
\lb{phi}
&&
\varphi_{(t_1,t_2,t_3)} \, :=\, \exp\bigl\{-{\frac{i \pi}{ 2\tau} \sum_{\alpha =1}^n
(t_1 x^2_\alpha  -t_2 y^2_\alpha+ t_3 z^2_\alpha)}\bigr\}.
\ea
Automorphisms $\phi_L$, $\phi_K$ and $\phi_F$ are elements of this family:
\be
\lb{phi-LK}
\phi_L\, =\, \phi_{(1,0,0)}\, , \qquad \phi_K\, =\, \phi_{(0,1,0)}\, , \qquad \phi_F\, =\, \phi_{(0,0,1)}\, .
\ee
\end{theor}

\ni
{\bf Proof.}
For the proof one uses decompositions of matrix units (\ref{projL}), (\ref{projK}), (\ref{projF}) and permutation formulas like ones in
(\ref{T-x}), (\ref{Omega-x}). An idea of the proof was suggested in \cite{IP.09}, Section 4.
\hfill$\blacksquare$
\smallskip

In view of remark \ref{rem3.3} and eq.(\ref{weaker}) one has

\begin{prop}
\lb{prop3.7}
Bi-invariant elements of the calculi are invariant under two-para\-met\-ric family of automorphisms $\phi_{(t_1,t_2,0)}$. R-traces of matrix monomials in matrices  $\Omega$ and $F$ are invariant under whole family of automorphisms $\phi_{(t_1,t_2,t_3)}$.
\end{prop}

\section{Gauss decomposition}
\lb{sec4}

We need one more structure to construct the unitary calculus, namely, the Gauss decomposition for the Lie derivatives.
To our knowledge such a decomposition is only known for the
RE algebras associated with the Drinfeld-Jimbo's R-matrix. So, from now on we consider the calculus associated with the R-matrix (\ref{DJ}).

Following \cite{FRT} we introduce two pairs of the RTT algebras generated by the upper/lower triangular matrices
 $L^{(+/-)}=|| \ell^{(\pm)}_{ij}||_{i,j=1}^{\;\;\;\;\;\;\, n}$, and
$K^{(+/-)}=|| k^{(\pm)}_{ij}||_{i,j=1}^{\;\;\;\;\;\;\, n}$ subject to the permutation relations
\ba
\lb{forL}
R L^{(\pm)}_2 L^{(\pm)}_1 \, =\, L^{(\pm)}_2 L^{(\pm)}_1  R, \;\,&\; &\;\,
R L^{(+)}_2 L^{(-)}_1 \, =\, L^{(-)}_2 L^{(+)}_1 R,
\\[2pt]
\lb{forK}
R K^{(\pm)}_2 K^{(\pm)}_1 \, =\, K^{(\pm)}_2 K^{(\pm)}_1  R, &\; &
R K^{(+)}_2 K^{(-)}_1 \, =\, K^{(-)}_2 K^{(+)}_1 R,
\ea
and to the $SL_q(n)$ reduction conditions\footnote{Here the quantum determinant for matrices $L^{(\mp)}$, $K^{(\mp)}$
is defined by formula (\ref{detT}), which is universal for the RTT algebras.}
\be
\lb{SL-forLK}
\Rdet L^{(\pm)}\, =\, \prod_{i=1}^n \ell_{ii}^{(\pm)}\, =\, 1, \qquad\quad
\Rdet K^{(\pm)}\, =\, \prod_{i=1}^n k_{ii}^{(\pm)}\, =\, 1,
\ee
and
\be
\lb{SL-forLK-2}
\ell^{(-)}_{ii} \ell^{(+)}_{ii}\, =\, k^{(-)}_{ii} k^{(+)}_{ii}\, =\, 1\quad \forall\; i=1,2,\dots ,n.
\ee

As is well known the RE algebra can be realized in terms of these upper/lower triangular RTT algebras (see, e.g., \cite{KSch}, pp.345-347).
So we do for the Lie derivatives $L$ and $K$:
\be
\lb{12}
L\, =\, q^{n-1/n} \bigl(L^{(-)}\bigr)^{-1} L^{(+)} \, , \qquad K\, =\, q^{n-1/n} K^{(+)} \bigl(K^{(-)}\bigr)^{-1} \, .
\ee
Note that normalization factor $q^{n-1/n}$ in these formulas
is necessary for compatibility of the $SL_q(n)$ reductions (\ref{L-sl}), (\ref{detK-sl}) and (\ref{SL-forLK}).
Indeed, one can calculate
$$
\Rdet L = q^{-1}\bigl(\Rdet L^{(-)}\bigr)^{-1} \Rdet L^{(+)},
$$
and the same for matrix $K$.\smallskip

An extension of the Gauss decomposition to the spectral variables is obvously central.
Less trivial is the extension for the algebras $\widehat{\goth D\goth C}_{gl/sl}[R]$.
It was elaborated in \cite{AF.92,STS,SWZ.93}. Below we present a list of permutation relations for the matrices $L^{(\pm)}$ and $K^{(\pm)}$
derived in  these papers
\ba
\lb{Lpm-T}
L^{(\pm)}_1 R^{\pm 1} T_1  &= & \gamma^{\pm 1} T_2\, P\, L^{(\pm)}_2 ,
\\
\lb{Kpm-T}
T_2 R^{\pm 1} K^{(\pm)}_2  &= & \gamma^{\pm 1} K^{(\pm)}_1\, P\, T_1 .
\ea
Recall that $\gamma = q^{1/n}$  and $P\in {\rm Aut}(V^{\otimes 2})$ is the permutation matrix.
\ba
\lb{L-K}
\bigl[ L^{(\pm)}_1  ,  K^{(\pm)}_2\bigr] &=&\bigl[ L^{(\pm)}_1  ,  K^{(\mp)}_2\bigr] \, =\, 0,
\\[2pt]
\lb{munurho-LK}
\bigl[ \,\xi  ,  L^{(\pm)}\bigr]  &=& \bigl[ \,\xi  ,  K^{(\pm)}\bigr]  \, =\, 0, \quad\forall \xi\in\{\mu_\alpha,\nu_\alpha,\rho_\alpha\},
\ea
\ba
\lb{Lpm-X}
L^{(\pm)}_1 R^{\pm 1} X_1  \,= \, X_2\,  L^{(\pm)}_1 R^{\pm 1}, && K^{(\pm)}_1\, X_2\, =\, X_2\, K^{(\pm)}_1
\quad \forall\, X= L,\Omega, F,
\\[1pt]
\lb{Kpm-Y}
Y_2 R^{\pm 1} K^{(\pm)}_2  \, = \, R^{\pm 1} K^{(\pm)}_2\, Y_1 , &&\;\; L^{(\pm)}_1\, Y_2\, =\, Y_2\, L^{(\pm)}_1
\quad\;\; \forall\, Y= K,\Theta .
\ea
One can  check that these relations are
(i) consistent with the previously defined permutation relations for $L$ and $K$, and (ii) respect reduction conditions
(\ref{SL-forLK}), (\ref{SL-forLK-2}).

\section{Unitary anti-involution}
\lb{sec5}

Now we are are ready to construct a unitary  anti-involution on $\widehat{\goth D\goth C}_{gl}[R]$.
\smallskip

We fix value of the quantization parameter $q$ on a unit circle: $q = e^{2 \pi i \tau}$, $\tau\in {\Bbb R}$.
In this case the Hermite conjugate of the Drinfeld-Jimbo R-matrix is
\be
\lb{Rdag}
R^{\dagger} = P R^{-1} P.
\ee

As a starting point we take the Hermite conjugation of the triangular components of the Lie derivatives
adopted in \cite{AF.92}:
\be
\lb{21}
(L^{(\pm)})^{\dagger} = (L^{(\mp)})^{-1}\, , \qquad
(K^{(\pm)})^{\dagger} = (K^{(\mp)})^{-1} \, ,
\ee
where by  "$\dagger$" we understand composition of the anti-linear algebra anti-involution and the matrix transposition.
It is easy to check that this setting is compatible with permutations relations and reduction conditions (\ref{forL})-(\ref{SL-forLK-2})
for $L^{(\pm)}$ and $K^{(\pm)}$.

\begin{rem}
The RTT algebras generated by matrices $L^{(\pm)}$ and $K^{(\pm)}$ can be endowed with the standard Hopf structure (\ref{bi-RTT})-(\ref{T-inverse}).
However, the $\dagger$ structure (\ref{21}) does not make them exactly the Hopf *-algebras. Instead, the compatibility condition for the coproduct and the
Hermite conjugation reads
$$
\left(\Delta X\right)^{\dagger\otimes\dagger}\, =\, \sigma\circ \Delta(X^\dagger), \qquad \forall X= L^{(\pm)}, K^{(\pm)},
$$
where $\sigma$ is the transposition map (see, e.g., \cite{Maj}, p.101).
\end{rem}

\subsection{Conjugation of spectral variables  $\mu_\alpha$, $\nu_\alpha$}\hfill  \par\vspace{3pt}
\lb{subsec5.1}

\noindent
In this subsection we calculate an effect of
Hermitean conjugation on RE algebras of Lie derivatives and on their spectra.

For the matrices of generators their Hermite conjugates look like
\be
\lb{LK-Hermite}
L^\dagger \, =\, L^{(\pm)} L^{-1} \left(L^{(\pm)}\right)^{-1},\qquad
K^\dagger \, =\, \left(K^{(\pm)}\right)^{-1} K^{-1} K^{(\pm)}\, .
\ee

Consider a set of bi-invariant elements, called {\em power sums}
\be
\lb{power-sums}
p_k := \Rtr L^k, \quad p_{-k} :=  \Rinvtr L^{-k}\, =\, q^{2n}\, \Rtr L^{-k}\quad \mbox{\em (see (\ref{n2}))}.
\ee
They are related with the coefficients of characteristic polynomial (\ref{charL}) by a $q$-version of Newton relations \cite{GPS.97}
\be
\lb{Newton}
p_k\, -\, q a_1 p_{k-1}\, +\dots +\, (-q)^{k-1} a_{k-1} p_1 \, +\, (-1)^k k_q a_k\, =\, 0, \qquad \forall \, k\geq 1.
\ee
These formulas are helpful for the following
\begin{prop}
On elements $p_k$, $a_k$, $b_k$ conjugation gives
\be
\lb{pa-Hermite}
p_k^\dagger\, =\, p_{-k}, \qquad a_k^\dagger\, =\,  a_{n-k}/a_n,\qquad b_k^\dagger\, =\,  b_{n-k}/b_n,
\ee
\end{prop}
\smallskip

\ni
{\bf Proof.}
Formulas for the power sums are obtained by a direct calculation:
\ba
\nn
p_k^\dagger &=& q^{2n}\, \Roptr ( L^\dagger)^k\, =\, q^{2n} \,\RopTr{1} \RTr{2} \,\underline{L^{(+)}_1 R L_1^{-k}} ( L_1^{(+)})^{-1}
\\
\nn
&=& q^{2n}\, \RopTr{1} \RTr{2}\, L_2^{-k} \underline{L_1^{(+)} R (L_1^{(+)})^{-1}}\, =\, q^{2n}\, \Rtr L^{-k}\, =\ p_{-k}.
\ea
Here we used permutation relations and formulas (\ref{n1}) from the Appendix.
For clarity we underlined expressions which are transformed in the next step.

To get conjugation formulas for $a_k$ we use the characteristic identity (\ref{charL}). Multiplying it by $L^{-k}$ and taking
$\Rinvtr\! =q^{2n} \Rtr$ we obtain
\ba
\nn
&
q^{2n}\left( p_{n-k} - q a_1 p_{n-k-1}+\dots + (-q)^{n-k-1} a_{n-k-1} p_1 \right)\, +\, (-q)^{n-k} a_{n-k} \Rinvtr I
&
\\[1pt]
\nn
&
+\, (-q)^{n-k+1} a_{n-k+1} p_{-1} +\dots + (-q)^n a_n p_{-k}\, =\, 0.
&
\ea
Simplifying the first term in brackets with the help of (\ref{Newton}) and using (\ref{n1}) for the second term
we find after collecting similar terms
\ba
\nn
\textstyle
p_{-k} -\frac{1}{q} \frac{a_{n-1}}{a_n} p_{-k+1} +\dots +(-\frac{1}{q})^{k-1}\frac{a_{n-k+1}}{a_n} p_{-1} + (-1)^k k_q \frac{a_{n-k}}{a_n}\, =\, 0.
\ea
Comparing this formula with the result of the conjugation of (\ref{Newton}) we conclude $a_k^\dagger = a_{n-k}/a_n$.

Formulas for $b_k^\dagger$ are obtained in the same way.
\hfill$\blacksquare$
\smallskip

The proposition together with the parameterization formulae (\ref{aL-elementary}), (\ref{bK}) suggests
\begin{cor}
\lb{cor}
In the spectrally extended RE algebras conjugation rules (\ref{LK-Hermite})
can be consistently complemented by
\be
\lb{munu-Hermite}
\mu_\alpha^\dagger\, =\, \mu_\alpha^{-1}, \qquad
\nu_\alpha^\dagger\, =\, \nu_\alpha^{-1}.
\ee
Hermite conjugation of the corresponding matrix idempotents reads
\be
\lb{PQ-Hermite}
(P^\alpha)^\dagger\, =\, L^{(\pm)} P^\alpha\! \left(L^{(\pm)}\right)^{-1},
\qquad
(Q^\alpha)^\dagger\, =\, \left(K^{(\pm)}\right)^{-1}\! Q^\alpha K^{(\pm)}.
\ee
\end{cor}
\smallskip

\subsection{Conjugation ansatz for $T$, $F$ and $\rho_\alpha$}\hfill  \par\vspace{3pt}
\lb{subsec5.2}

\noindent
Formula for Hermite conjugation of $T$ in the Heisenberg double algebra was suggested in \cite{AF.92}. Generalizing it for the
differential calculi algebras $\widehat{\goth D\goth C}_{gl/sl}[R]$ we write down following ansatz
\be
\lb{ansatzT}
T^\dagger\, =\, q^{n-1/n} \bigl(K^{(-)}\bigr)^{-1} \widetilde{T^{-1}} \bigl(L^{(-)}\bigr)^{-1}.
\ee
Here we use  shorthand notation $\widetilde{X}$ for the image of $X$ under some automorphism from the family (\ref{phi}).
Its explicit form is to be specified later on. Our choice of numeric factor $q^{n-1/n}$  will also be argued below.

The suggested $T^\dagger$ has to satisfy the Hermite conjugates of permutation relations (\ref{RTT}), (\ref{Lpm-T}), (\ref{Kpm-T}) 
\ba
\lb{RTTdag}
R\, T_1^\dagger T_2^\dagger & =& T_1^\dagger T_2^\dagger R ,
\\
\lb{Lpm-Tdag}
T_2^\dagger R^{\pm 1} L^{(\pm)}_2  &= & \gamma^{\pm 1} L^{(\pm)}_1\, P\, T_1^\dagger ,
\\
\lb{Kpm-Tdag}
K^{(\pm)}_1 R^{\pm 1} T_1^\dagger  &= & \gamma^{\pm 1} T_2^\dagger\, P\, K^{(\pm)}_2 ,
\ea
It is a standard exercise in R-matrix calculations to verify these equalities. We only mention that while proving (\ref{RTTdag}) one finds  a remarkable relation~
\be
\lb{interest1}
T^\dagger_1\, \widetilde{T_2} = \widetilde{T_2}\, T^\dagger_1.
\ee

It is also straightforward to test compatibility of the ansatz with  permutation relations (\ref{T-mu}), (\ref{T-nu}) of the spectrally extended algebras
$\widehat{\goth D\goth C}_{gl/sl}[R]$.

Less easy is checking consistency of the ansatz with $SL_q(n)$ reduction condition (\ref{sl-1}). It is suitable to consider the Hermite conjugate of its inverse
\be
\lb{RdetTdag}
\bigl( \Rdet T^{-1} \bigr)^\dagger = 
\Tr{1, \dots ,n} \bigl( A^{(n)} (T^\dagger_n)^{-1}\dots  (T^\dagger_1)^{-1}\bigr),\quad \mbox{\em (see (\ref{n5})).}
\ee
To calculate it we separate factors $K^{(-)}$, $\widetilde{T}$ and $L^{(-)}$ in the expression
\ba
\nn
&
(T^\dagger_n)^{-1}\dots  (T^\dagger_1)^{-1}\, =\,
\eta^{-n}(L^{(-)}_n \widetilde{T_n} K^{(-)}_n) \dots (L^{(-)}_1 \widetilde{T_1} K^{(-)}_1)
&
\\[2pt]
\nn
&
=\, \eta^{-n} \gamma^{n(n-1)} \left(\prod_{i=1}^n J_i\right)^{-1}\bigl(L^{(-)}_n\dots  L^{(-)}_1\bigr) \bigl(\widetilde{T_1}\dots  \widetilde{T_n}\bigr)
\bigl(K^{(-)}_n\dots  K^{(-)}_1\bigr), 
&
\ea
where by $\eta$ we denote the numeric factor in the ansatz: $\eta = q^{n-1/n}$.
Substituting this expression in (\ref{RdetTdag}), evaluating  Jucys-Murphy elements $J_i$ on $A^{(n)}$ and using the rank$=1$
property of $A^{(n)}$ we obtain
$$
\bigl( \Rdet T^{-1} \bigr)^\dagger = 
\bigl(\eta^{-1} \gamma^{n-1} q^{n-1}\bigr)^n\, \Rdet L^{(-)}\, \Rdet \widetilde{T}\, \Rdet K^{(-)}\,  =\,1,
$$
where the last equality is satisfied due to $SL_q(n)$ reduction conditions and due to our choice of normalization $\eta$ in the ansatz.
\smallskip

Now we discuss Hermite conjugation of matrix $F$,
postponing investigation of involutivity of the ansatz (\ref{ansatzT}) to subsection \ref{subsec5.5}.

Given formulae (\ref{LK-Hermite}) for $L^\dagger$ and $K^\dagger$ and the ansatz for $T^\dagger$ one can calculate $F^\dagger$
\be
\lb{Fdag}
F^\dagger\, =\, L^{(-)} \widetilde{F^{-1}} \bigl( L^{(-)}\bigr)^{-1}.
\ee
This expression is quite similar to those for $L^\dagger$, $K^\dagger$ and hence,
considerations of section \ref{subsec5.1}  can be repeated with little modifications for the eigenvalues of $F$.
We collect the results in 

\begin{prop}
\lb{prop5.4}
Formula (\ref{Fdag}) for Hermite conjugation of matrix $F$  
determines conjugation rules for the coefficients of its characteristic polynomial
\be
\lb{c-dag}
c_k^\dagger\, =\, c_{n-k}/c_n\, .
\ee
These rules in turn, agree with the unitary conjugation prescriptions for $F$'s eigenvalues
\be
\lb{rho-dag}
\rho_\alpha^\dagger \, =\, \rho_\alpha^{-1},
\ee
which result in following Hermite conjugation for their corresponding matrix idempotents
\be
\lb{S-dag}
(S^\alpha)^\dagger\, =\, L^{(-)} \widetilde{S^\alpha}\! \left(L^{(-)}\right)^{-1}.
\ee
\end{prop}
\smallskip

\ni
{\bf Proof.}
The only point which needs to be commented here is invariance of the r.h.s. of (\ref{c-dag}) under the automorphism
from the ansatz. This fact follows by proposition \ref{prop3.7}.  
\hfill$\blacksquare$
\smallskip

\subsection{Conjugation ansatz for differential forms}\hfill  \par\vspace{3pt}
\lb{subsec5.4}

\noindent
In this subsection we introduce an ansatz for Hermite congutation in the algebra of differential forms. 
From now on things start to be different in cases  $GL_q(n)$ and $SL_q(n)$.
We show briefly consistency of the ansatz with the algebra structure of $\widehat{\goth D\goth C}_{gl}[R]$ and consider possibility
of the $SL_q(n)$ reduction of the conjugation. 

For Hermite conjugate of the matrix $\Omega^g$ we write down the following ansatz
\be
\lb{ansatzO}
(\Omega^g\!)^\dagger\, =\, -\, L^{(-)} \widetilde{F^{-1}\Omega^g}\, \bigl( L^{(-)}\bigr)^{-1}\, .
\ee
Here by $\widetilde{\dots}$ we denote an action of the same automorphism as in (\ref{ansatzT}).

It is not hard to verify that the matrix elements of ${(\Omega^g)}^\dagger$ indeed satisfy the conjugated permutation relations
(\ref{Lpm-X}), (\ref{OmegaT}), (\ref{OmegaOmega-GL}):
\ba
\lb{LKpm-Odag}
{\Omega^g_2}^\dagger R^{\pm 1} L^{(\pm)}_2  \, = \, R^{\pm 1} L^{(\pm)}_2\, {\Omega^g_1}^\dagger , && K^{(\pm)}_1\, {\Omega^g_2}^\dagger\, =\,
{\Omega^g_2}^\dagger\, K^{(\pm)}_1,
\\
\lb{Odag-Tdag}
T_2^\dagger R\, {\Omega_2^g}^\dagger R^{-1}& =& {\Omega_1^g}^\dagger\, T_2^\dagger\,   ,
\\
\lb{O-Odag}
R\,{\Omega^g_2}^\dagger R^{-1} {\Omega^g_2}^\dagger &=& -\, {\Omega^g_2}^\dagger R^{-1} {\Omega^g_2}^\dagger R^{-1} \, .
\ea
Here as an intermediate step in proving eq.(\ref{Odag-Tdag}) one obtains a remarkable commutativity relation
\be
\lb{interest2}
{\Omega^g_1}^\dagger\, \widetilde{T_2}\, =\, \widetilde{T_2}\, {\Omega_1^g}^\dagger .
\ee
Permutation relations for spectral variables $\mu_\alpha$, $\nu_\alpha$, $\rho_\alpha$ (\ref{Omega-mu})-(\ref{weaker})
are also compatible with the ansatz. 

Consider now action of the conjugation on the subalgebra generated by the R-traceless forms (\ref{O-SL}).
In view of (\ref{n5}), under conjugation they go into the $R_{op}$-traceless matrices
\be
\lb{Odag-sl}
\Omega^\dagger\, =\, (\Omega^g)^\dagger -\frac{q^n}{n_q} \Roptr\! \bigl(\Omega^g\bigr)^\dagger I,
\ee
which satisfy permutation relations
\ba
\lb{O-Odag-sl}
R^{-1}\Omega^\dagger_2 R^{-1} \Omega^\dagger_2 \, +\, \Omega^\dagger_2 R^{-1} \Omega^\dagger_2 R & =&
\kappa_{1/q} \bigl( (\Omega^\dagger_2)^2\, +\, R^{-1}(\Omega^\dagger_2)^2 R^{-1}\bigr)\, .
\ea
and hence,  generate closed subalgebra in the external algebra (\ref{OmegaOmega-GL}), as well as $\Omega$ did.
However, the subalgebra generated by $\Omega^\dagger$ goes beyond the $SL_q(n)$ differential calculus  described earlier. 
Indeed, with the use of R-techniques one can express $R_{op}$-traces of the conjugated 1-forms in terms of the R-traces of non-conjugated ones  
\be
\lb{tr-Odag}
\Roptr\! \bigl(\Omega^g\bigr)^\dagger\, =\, -\, \Rtr \bigl( F^{-1} \Omega^g\bigr)
\ee
Namely an appearence of the matrix factor $F^{-1}$ in this formula shows clearly that the conjugation map $\dagger$ does not preserve
the $SL_q(n)$ differential calculus. So we are left in a situation  where two different mutually conjugate $SL_q(n)$ calculi subalgebras
lie inside the $GL_q(n)$ calculus algebra.
\smallskip

\subsection{Involutivity}\hfill  \par\vspace{3pt}
\lb{subsec5.5}

\noindent
In this subsection we fix uniquely the automorphism $\widetilde{\dots}$ in the ansatz (\ref{ansatzT}), (\ref{ansatzO}) demanding involutivity of
the conjugation $\dagger$. We then summarize considerations of the present section in a theorem.\smallskip

Here we expand notation 
\be
\lb{tilde-not}
\widetilde{X}\, :=\, \varphi X \varphi^{-1},
\ee
where $\varphi$ is one of automorphism's generating elements (\ref{phi}).
We assume $\varphi^\dagger = \varphi^{-1}$ keeping in mind unitarity of the spectral variables. We now calculate $T^{\dagger\dagger}$ and 
$\Omega^{\dagger\dagger}$:

\ba
\nn
\bigl(T^\dagger\bigr)^\dagger &=& q^{1/n-n}\, \bigl(L^{(+)}\bigr)^{-1}\, (\varphi^{-1})^\dagger\, \bigl(T^{-1}\bigr)^\dagger\,
\varphi^\dagger\, \bigl(K^{(+)}\bigr)^{-1}
\\
\nn
&=& q^{1/n-n}\,\bigl(L^{(+)}\bigr)^{-1}\, \varphi\, \bigl( q^{1/n-n}\, L^{(-)}\, \varphi\, T\, \varphi^{-1}\, K^{(-)}\bigr) \varphi^{-1}\,
\bigl(K^{(+)}\bigr)^{-1}
\\
\nn
&=& L^{-1}\, \varphi^2\, T\, \varphi^{-2}\, K^{-1}\, ;
\\[5pt]
\nn
\bigl({\Omega^g}^\dagger\bigr)^\dagger &=& \bigl(L^{(+)}\bigr)^{-1}\, (\varphi^{-1})^\dagger\, \bigl(\Omega^g\bigr)^\dagger
\bigl(F^{-1}\bigr)^\dagger\, \varphi^\dagger\, L^{(+)}
\\
\nn
&=&
\bigl(L^{(+)}\bigr)^{-1}\, \varphi\, \bigl( L^{(-)}\, \varphi\, F^{-1}\, \Omega^g \varphi^{-1} (L^{(-)})^{-1} \bigr)
\bigl( L^{(-)} \varphi F \varphi^{-1} (L^{(-)})^{-1}\bigr)\, \varphi^{-1}\, L^{(+)}
\\
\nn
&=&
L^{-1} F^{-1} \varphi^2\, \Omega^g\, \varphi^{-2}\, F L.
\ea
So we conclude that conditions $T^{\dagger\dagger}=T$ and $(\Omega^g)^{\dagger\dagger}=\Omega^g$ are satisfied with the choice
\be
\lb{ansatz-phi}
\varphi^{2}\, =\, \varphi_{(1,1,1)}, \qquad \mbox{ that is}\qquad
\varphi\, =\, \exp\bigl( -\frac{i \pi}{4 \tau} \sum_{\alpha=1}^n (x^2_\alpha - y^2_\alpha +z^2_\alpha)\bigr).
\ee
Now we are ready to formulate final
\begin{theor}
\lb{theor5.5}
For the Drinfeld-Jimbo R-matrix (\ref{DJ}) consider spectrally extended algebra 
$\widehat{\goth D\goth C}_{gl}[R]$ of the differential calculus over $GL_q(n)$
taking parameter $q$ on a unit circle: $q = e^{2 \pi i \tau}$, $\tau \in {\Bbb R}$.

The anti-linear algebra anti-homomorphism given on the generators by formulas
\begin{center}
(\ref{21}), (\ref{munu-Hermite}), (\ref{ansatzT}), (\ref{Fdag}), (\ref{rho-dag}), (\ref{ansatzO}), (\ref{tilde-not}), (\ref{ansatz-phi})
\end{center}
defines unitary type anti-involution on $\widehat{\goth D\goth C}_{gl}[R]$. This unitary structure respects the bicovariance
property of the calculus in a sense that the algebra $\widehat{\goth D\goth C}_{gl}[R]$ can be endowed with the structure of
the bicovariant bimodule over Hopf algebra ${\cal F}^\dagger[R]$ generated by the matrix components of $T^\dagger$. The left and right 
${\cal F}^\dagger[R]$ coactions $\delta^\dagger_{\ell/r}$ are defined on the generators as
\be
\lb{coacTdag}
\delta^\dagger_{\ell/r} (T^\dagger_{ij})\, =\, {\textstyle\sum_{k=1}^n}T^\dagger_{ik}\otimes T^\dagger_{kj},
\ee
\ba
\lb{coacXdag}
\delta^\dagger_{r} (X^\dagger_{ij})&=& {\textstyle\sum_{k,p=1}^n} \bigl(1\otimes (T^{-1})^\dagger_{ik}\bigr)\bigl(X^\dagger_{kp}\otimes 1\bigr)
\bigl(1\otimes T^\dagger_{pj}\bigr),
\quad \delta^\dagger_{\ell} (X^\dagger_{ij}) = X^\dagger_{ij}\otimes 1,\qquad
\\
\lb{coacYdag}
\delta^\dagger_{\ell} (Y^\dagger_{ij})&=& {\textstyle\sum_{k,p=1}^n} \bigl(T^\dagger_{ik}\otimes 1\bigr) \bigl(1\otimes Y^\dagger_{kp}\bigr)
\bigl((T^{-1})^\dagger_{pj}\otimes 1\bigr), \quad \;\,
\delta^\dagger_{r} (Y^\dagger_{ij})\, =1\otimes Y^\dagger_{ij},\qquad
\ea
where  $X^\dagger = {\Omega^g}^\dagger, L^\dagger, F^\dagger, \dots$; $Y^\dagger = {\Theta^g}^\dagger, K^\dagger, \dots$ . Naturally, one can  consider
${\cal F}^\dagger[R]$ coactions $\delta^\dagger_{\ell/r}$ as Hermite conjugates of the ${\cal F}[R]$ coactions $\delta_{r/\ell}$, respectively.

Restriction of conjugation $\dagger$ to the subalgebra $\widehat{\goth D\goth C}_{sl}[R]$ results in the involutive anti-homomorphism
of the two $SL_q(n)$ type subalgebras, generated by the $R/R_{op}$-traceless matrices of 1-forms $\Omega$ (\ref{O-SL}) and $\Omega^\dagger$
(\ref{Odag-sl}), respectively.
\end{theor}
\smallskip

{\bf Acknowledegment}  \par
I thank  Alexei Isaev, Ludwig Faddeev, Oleg Ogievetsky, Dimitry Gurevich  and Pavel Saponov for inspirating and helpful discussions,
sharing ideas and for the years of fruitful collaboration.

\appendix

\numberwithin{equation}{section}

\section{R-matrices}
\label{append}

Throughout the paper we consider various matrices acting on tensor powers of some finite dimensional
vector space $V$. For these matrices we use, by now standard, compressed matrix notation.
Namely, with any matrix $X\in {\rm End} (V^{\otimes k})$ we associate series of matrices $X_i\in {\rm End} (V^{\otimes n})$, $n\geq k$,
\be
\lb{X-i}
X_i\, := I^{\otimes (i-1)}\otimes X\otimes I^{\otimes (n-k+1-i)} ,\qquad i=1,2,\dots , n-k+1,
\ee
where $I\in {\rm Aut}(V)$ is the identity. For $X\in {\rm End}(V^{\otimes 2})$ in certain occasions we also use notation $X_{ij}$ for the matrices
acting nontrivially in spaces with labels $i$ and $j$, $i<j$. In these notation $X_i \equiv X_{i,i+1}$.\smallskip

An operator $R\in {\rm Aut}(V^{\otimes 2})$ satisfying braid relation
\be
\lb{ybe}
R_1 R_2 R_1\, =\, R_2 R_1 R_2\, ,
\ee
is called an {\em {\rm R}-matrix}.
Permutation $P$: $P(u\otimes v) = v\otimes u$,  is the R-matrix.
If $R$ is the R-matrix, so are the operators $R^{-1}$ and  $R_{op} := P R P$.\smallskip

An {\rm R}-matrix  $R$ is called {\em skew invertible} if there exist an
operator $\Rpsi \in {\rm End}(V^{\otimes 2})$  such that
\be
\lb{psi}
\Tr{2} R_{12} \Rpsi_{23} = \Tr{2} \Rpsi_{12} R_{23} = P_{13} \; .
\ee
Here $\tr_{\!(i)}$ denotes trace operation in $i$-th space.
With a skew invertible R-matrix $R$ one associates matrix
$\D \in {\rm End}(V)$:~
$${\D}_1 := \Tr{2} \Rpsi_{12},$$
by which one defines a notion of {\em {\rm R}-trace}, $\Rtr$. Namely, for any $X\in {\rm End}(V)$
$$
\Rtr(X) := \tr(\D X) .
$$
The operation $\Rtr$ is often called a quantum trace or, shortly, a $q$-trace.
We use the name R-trace to emphasize dependence of this operation on a choice of the R-matrix.
Properties of the R-trace are listed in \cite{IP.09}, Sec. 2.2.

An R-matrix $R$ whose minimal polynomial is quadratic
is called {\em Hecke type}. By an appropriate rescaling
one can turn its minimal polynomial to a form
\be
\lb{H}
(R -q I)(R+q^{-1}I) = 0\, ,
\ee
known under the name {\em Hecke condition}. Skew invertible Hecke type R-matrices are used for quantizing differential geometric constructions over linear (super)groups.\smallskip

To specify $GL/SL(n)$ cases we impose conditions on the R-matrix.
First, we demand that  the R-matrix eigenvalue $q\in {\Bbb C}\setminus \{0\}$ does not coincide with certain roots of unity:
\be
\lb{restrict}
i_q := {(q^i-q^{-i})/(q-q^{-1})}\neq 0 \quad \forall i=2,3,\dots ,n.
\ee
In this case by the Hecke type R-matrix one can construct series of idempotents $A^{(i)}\in {\rm End}(V^{\otimes i})$, $i=1,\dots ,n$,
called {\em $q$-antisymmetrizers}. Their inductive definition reads
\be
\lb{q-anti}
A^{(1)} = I, \qquad
A^{(i)} = \frac{(i-1)_q}{i_q}\,
A^{(i-1)}\,
\Bigl(\frac{q^{i-1}}{(i-1)_q}\, I\, -\, R_{i-1}\Bigr)
A^{(i-1)},
\ee
and their properties are listed in \cite{IP.09}, Sec. 2.4.

A skew invertible Hecke type R-matrices whose eigenvalues satisfy (\ref{restrict}) is called {\em $GL(n)$ type}
if conditions
\ba
\lb{GL(n)}
A^{(n)} \Bigl( \frac{q^n}{ n_q}\,I\,-\,R_n\Bigr)A^{(n)}&=& 0,\hspace{57mm}
\\[-2mm]
\nn
\mbox{and}\hspace{77mm}&&
\\[-2mm]
\lb{GL(n)-2}
\qquad {\rm rk}A^{(n)} &=& 1
\ea
are fulfilled.
The R-matrix is called {\em $SL(n)$ type} if additionally condition
\be
\lb{SL(n)}
\Tr{2,\dots ,n+\! 1}\left(P_1 P_2\dots P_n\, A^{(n)}\right)\, \propto\, I_1
\ee
is satisfied. The latter condition guarantees centrality of the element $\Rdet T$ in the differential calculus algebra
${\goth D\goth C}[R]$ \cite{Gur} (see also \cite{IP.09}) and thus, makes the $SL_q(n)$ reduction possible.\smallskip

If the R-matrix $R$ is $GL(n)$/$SL(n)$ type, then so are the R-matrices $R^{-1}$ and $R_{op}$.\smallskip

We complete the Appendix with the list of formulas which are valid for the R-traces of the $GL(n)$ type R-matrices.
\ba
\lb{n1}
&
\RTr{2} R_{12} = I_1, \quad\; {\rm Tr}_{\hspace{-2pt} \raisebox{-2pt}{\tiny\sl $R_{op}$}\,\raisebox{2pt}{\scriptsize$(1)$}}
 R_{12}  = I_2, \quad\; \Rtr I = \Roptr I = q^{-n} n_q,
&
\\[2pt]
\lb{n2}
&
D_{R_{op}}\, =\, q^{-2n} (\D)^{-1}, \qquad D_{R^{-1}}\, =\, q^{2n} \D,
&
\ea
\be
\lb{n3}
\Tr{1,\dots n}\left(A^{(n)}\dots \right) = q^{n^2}\,\RTr{1,\dots n}\left(A^{(n)}\dots \right) =
q^{n^2}\,{\rm Tr}_{\hspace{-2pt} \raisebox{-2pt}{\tiny\sl $R_{op}$}\,\raisebox{2pt}{\scriptsize$(1, \dots ,n)$}}
\left(A^{(n)}\dots \right)
\ee

For the Drinfeld-Jimbo R-matrix (\ref{DJ}) explicit expressions for the matrices of R-traces are
\ba
\lb{n4}
\D\, =\, \mbox{diag}\{q^{1-2n},q^{3-2n},\dots ,q^{-1}\}, \quad
\Dop\, =\, \mbox{diag}\{q^{-1},q^{-3},\dots ,q^{-2n+1}\},
\ea
and in case $|q|=1$ one has
\be
\lb{n5}
\left(\Rtr\right)^\dagger \, =\, q^{2n}\, \Roptr  , \qquad
\bigl( A^{(n)}\bigr)^\dagger \, =\, P^{(n)} A^{(n)} P^{(n)},
\ee
where $P^{(n)}:=P_1 (P_2 P_1)\cdot\dots\cdot (P_{n-1}\dots P_2 P_1)$ is the operator inversing enumeration of vector spaces in $V^{\otimes n}$.



\end{document}